\documentclass[12pt]{iopart}
\usepackage[utf8]{inputenc}
\usepackage[toc]{appendix}
\expandafter\let\csname equation*\endcsname\relax
\expandafter\let\csname endequation*\endcsname\relax
\usepackage{amsfonts}
\usepackage{bm}
\usepackage{cite}
\usepackage{empheq}
\usepackage{enumerate}
\usepackage{hyperref}
\hypersetup{
    colorlinks,
    linkcolor={blue!90!black},
    citecolor={blue!90!black},
    urlcolor={blue!80!black}
  }

\usepackage{float}
\usepackage{amsthm}

\usepackage{mathtools}
\usepackage{mathrsfs}
\usepackage[new]{old-arrows}
\usepackage{xcolor}
\usepackage{tikz}
\usetikzlibrary{positioning, arrows.meta}
\usepackage{amssymb,color}
\usepackage{cancel}
\usepackage[normalem]{ulem}

\theoremstyle{plain}
\newtheorem{theorem}{Theorem}[section]

\theoremstyle{definition}
\newtheorem{lemma}[theorem]{Lemma}
\newtheorem{remark}[theorem]{Remark}

\renewcommand{\star}{*}
\renewcommand{\langle}{(}
\renewcommand{\rangle}{)}

\newcommand{\bb}{\mathbb}
\newcommand{\mrm}{\mathrm}
\newcommand{\mc}{\mathcal}
\newcommand{\wt}{\widetilde}

\newcommand{\rd}{\mathrm{d}}

\newcommand{\nel}{{N_{\mathrm{el}}}}
\newcommand{\elindex}{k}
\newcommand{\elindextwo}{\ell}
\newcommand{\elindexthree}{{\elindextwo'}}
\newcommand{\elindexfour}{\elindex'}

\newcommand{\el}{x}

\newcommand{\nuc}{X}

\newcommand{\vel}{V}

\newcommand{\acc}{A}
\newcommand{\mom}{P}

\newcommand{\wf}{\psi}

\newcommand{\lap}{\Delta_\el}

\newcommand{\nab}{\nabla_\el}

\newcommand{\dens}{\rho}

\newcommand{\nnuc}{{N_{\mathrm{nuc}}}}
\newcommand{\nucindex}{K}
\newcommand{\nucindextwo}{L}

\newcommand{\mass}{M}
\newcommand{\charge}{Z}

\newcommand{\schop}{H}
\newcommand{\schoplin}{H^{\mathrm{lin}}}

\newcommand{\pot}{V}

\newcommand{\extpot}{\pot_{\mathrm{ext}}}

\newcommand{\potX}{\pot_{\mathrm{x}}}

\newcommand{\energy}{E}

\newcommand{\potH}{\pot_{\mathrm{H}}}

\newcommand{\potC}{\pot_{\mathrm{c}}}

\newcommand{\ltwoN}{L^2 (\bb{R}^3;\bb{C}^{\nel}  )}

\newcommand{\honeN}{H^1 (\bb{R}^3;\bb{C}^{\nel}  )}

\newcommand{\htwoN}{H^2 (\bb{R}^3;\bb{C}^{\nel}  )}

\newcommand{\KSCL}{$\mrm{\eqref{eq:KSN}}$}
\newcommand{\Bnuc}{\mathcal{B}_{\mrm{nuc}}(\tau) }
\newcommand{\Bel}{\mathcal{B}_{\mrm{el}}(\tau) }

\numberwithin{equation}{section}

\makeatletter
\newcommand{\subalign}[1]{%
  \vcenter{%
    \Let@ \restore@math@cr \default@tag
    \baselineskip\fontdimen10 \scriptfont\tw@
    \advance\baselineskip\fontdimen12 \scriptfont\tw@
    \lineskip\thr@@\fontdimen8 \scriptfont\thr@@
    \lineskiplimit\lineskip
    \ialign{\hfil$\m@th\scriptstyle##$&$\m@th\scriptstyle{}##$\hfil\crcr
      #1\crcr
    }%
  }%
}
\makeatother

\begin{document}
\title[Solutions to the TDKS equations coupled with classical nuclear dynamics]{Local existence and uniqueness of solutions to the time-dependent Kohn--Sham equations coupled with classical nuclear dynamics}
\author{Bj\"orn Baumeier$^{1,2}$, Onur \c{C}aylak$^{1,2}$, Carlo Mercuri$^3$, Mark Peletier$^{1,2}$, Georg Prokert$^1$ and Wouter Scharpach$^{1,2}$}
\address{$^1$ Department of Mathematics and Computer Science, Eindhoven University of Technology, P.O. Box 513, 5600 MB Eindhoven, The Netherlands}
\address{$^2$ Institute for Complex Molecular Systems, Eindhoven University of Technology, P.O. Box 513, 5600 MB Eindhoven, The Netherlands}
\address{$^3$ Dipartimento di Fisica, Informatica e Matematica, Universit\`a degli Studi di Modena e Reggio Emilia, Via Campi 213/b, Modena, Italy}
\ead{b.baumeier@tue.nl; o.caylak@tue.nl; carlo.mercuri@unimore.it; m.a.peletier@tue.nl; g.prokert@tue.nl; w.l.j.scharpach@tue.nl}
\begin{abstract}
We prove short-time existence and uniqueness of solutions to the initial-value problem associated with a class of time-dependent Kohn--Sham equations coupled with Newtonian nuclear dynamics, combining Yajima's theory for time-dependent Hamiltonians with Duhamel's principle, based on suitable Lipschitz estimates. We consider a pure power exchange term within a generalisation of the so-called Local Density Approximation (LDA), identifying a range of exponents for the existence and uniqueness of $H^2$ solutions to the Kohn--Sham equations.  
\end{abstract}
\tableofcontents
\section{Introduction}
\subsection{Main result}
In this paper we study the existence and uniqueness of solutions to the system
\begin{subequations}
\label{eq:KSN}
\begin{empheq}[]{alignat=2}
i\dot{\wf}_\elindex
&=-\tfrac{1}{2}\lap\wf_\elindex-\sum_{\nucindex=1}^\nnuc\frac{\charge_\nucindex}{|\cdot-{}\nuc_\nucindex|}\wf_\elindex
+\biggr(\frac{1}{|\,\cdot\,|}\star\dens\biggr)\wf_\elindex +
\lambda\dens^{q-1}\wf_\elindex,\label{eq:KSNa}\\
\ddot{\nuc}_\nucindex
&=
\frac{\charge_\nucindex}{\mass_\nucindex}\biggr[\int{}\frac{\el-\nuc_\nucindex}{|\el-\nuc_\nucindex|^3}\dens(\el)\rd{}\el+\sum_{\nucindextwo=1,\nucindextwo\neq \nucindex}^\nnuc\charge_\nucindextwo\frac{\nuc_\nucindex-\nuc_\nucindextwo}{|\nuc_\nucindex-\nuc_\nucindextwo|^3}\biggr],\label{eq:KSNb}\end{empheq}
\end{subequations}
where $\nel,\nnuc\in\bb{N},$ $\charge_\nucindex\in\bb{N},$ $\mass_\nucindex\geq 0,$ $\lambda\in\bb{R},$  $q>1$ are given, and $\elindex = 1,\ldots,\nel$, $\nucindex =1,\ldots,\nnuc$.

We will use the short-hand notation $\wf= (\wf_1,\ldots,\wf_\nel )$ and $\nuc= (\nuc_1,\ldots,\nuc_\nnuc )\in\bb{R}^{3\nnuc}$, with
\begin{align*}
\wf_\elindex=\wf_\elindex(\el,t),\qquad  \nuc_{\nucindex}= \nuc_{\nucindex}(t),\qquad \el\in\bb{R}^3,\qquad t\geq 0.
\end{align*}
In the above equations, for all $\wf(\cdot):  [0,\tau_{\textrm{max}} )  \rightarrow  H^2 (\bb{R}^3;\bb{C}^\nel )$ we set
\begin{align*}
\dens=\sum_{\elindex=1}^\nel |\wf_\elindex |^2.
\end{align*}
Moreover, for all $\nuc\in\bb{R}^{3\nnuc}$ and $\elindex= 1,\ldots,\nel$, we define
\begin{align}
 (\schop[\nuc,\dens]\wf )_\elindex:=-\tfrac{1}{2}\lap\wf_\elindex-\sum_{\nucindex=1}^{\nnuc}\frac{\charge_{\nucindex}}{ |\cdot-\nuc_\nucindex |}\wf_\elindex+ \biggr(\frac{1}{ |\,\cdot\, |}\star\dens \biggr)\wf_\elindex
+
\lambda\dens^{q-1}\wf_\elindex.
\label{eq:KSHam}
\end{align}

The dynamics of the elements~$\nuc(\cdot):[0,\tau_{\textrm{max}} )  \rightarrow \bb{R}^{3\nnuc}$ is driven by the acceleration function $\acc=\acc^1+\acc^2$, whose components are defined as
\begin{align}
\acc^1_\nucindex[\dens](\nuc)&:=\frac{\charge_{\nucindex}}{\mass_{\nucindex}}\int\frac{\el-\nuc_{\nucindex}}{|\el-\nuc_{\nucindex}|^3}\dens(\el)\rd \el,
\qquad \acc^2_\nucindex(\nuc)&:=\frac{\charge_{\nucindex}}{\mass_{\nucindex}}\sum_{\nucindextwo=1,\nucindextwo\neq \nucindex}^\nnuc \charge_{\nucindextwo}\frac{\nuc_{\nucindex}-\nuc_\nucindextwo}{|\nuc_{\nucindex}-\nuc_\nucindextwo|^3}.\label{eq:accel}
\end{align}
The main result of this paper is the following.
\begin{theorem}
\label{thm:shorttimeexistence}
Let $q\geq 7/2$ and $\lambda\in\bb{R}$. Further, let $\wf^0\in H^2 (\bb{R}^3;\bb{C}^\nel )$, $\vel^0\in\bb{R}^{3\nnuc}$ and $\nuc^0\in\bb{R}^{3\nnuc}$ be given, with $\nuc^0_{\nucindex}\neq \nuc^0_{\nucindextwo}$ for $\nucindex\neq \nucindextwo$.

Then, there exists $\tau>0$ such that the initial-value problem associated with the system \KSCL{} with $\wf(0)=\wf^0$, $\nuc(0)=\nuc^0$ and $\dot{\nuc}(0)=\vel^0$ has a unique solution $(\wf,\nuc)\in \mathcal{X}(\tau),$
where \begin{align*}
\mathcal{X}(\tau):=C^1 ([0,\tau];L^2 (\bb R^3;\bb{C}^\nel ) )\cap C^0 ([0,\tau];H^2 (\bb{R}^3;\bb{C}^\nel ) )\times C^2 ([0,\tau];\mathbb{R}^{3\nnuc} ).
\end{align*}
\end{theorem}
\subsection{Physical motivation}
Problems such as \eqref{eq:KSN} describe the nonadiabatic dynamics of molecular, spin-unpolarised systems involving an even number  $\nel\in 2\bb{N}$ of electrons and $\nnuc\in\bb{N}$ nuclei with masses $\mass_1,\ldots,\mass_{\nnuc}$ and charges $\charge_1,\ldots,\charge_{\nnuc}$. See e.g. \cite{Hohenberg1964InhomogeneousGas,Lieb1983DensitySystems,Anantharaman2009ExistenceChemistry,Perdew2003DensityCentury,Cohen2012ChallengesTheory,Kohn1965Self-ConsistentEffects,Ullrich2012Time-DependentApplications,LeBris2005FromJourney,Perdew1981Self-interactionSystems}, which form a sample of the extensive body of literature on both physical and mathematical aspects of the so-called Density-Functional Theory (DFT), which comprises the framework of the Time-Dependent Kohn--Sham (TDKS) equations, given in \eqref{eq:KSNa}. These equations, which using \eqref{eq:KSHam} can be written as
\begin{align}
\label{eq:KS}      i\dot{\wf}&=\schop[\nuc,\dens]\wf,
      \end{align}
describe the electronic evolution in terms of single-particle wave functions $\wf_\elindex$, known in the physical literature as the Kohn--Sham (KS) orbitals. The TDKS equations have been extensively considered as an approximation to the time-dependent Schr\"odinger equation, which reduces the electronic dynamics to a single-particle description based on the KS density function $\dens$. For convenience, we briefly recall the physical interpretation of each potential in the KS Hamiltonian $\schop$ from \eqref{eq:KSHam}, which can be written as
\begin{align}
\schop[\nuc,\dens]=-\tfrac{1}{2}\lap+\extpot[\nuc]+\pot_\text{Hxc}[\dens],\qquad \pot_{\text{Hxc}}:=\potH+\potX+\potC.\label{eq:Ham}
\end{align}
The different terms appearing in \eqref{eq:Ham} are defined as follows.
\noindent The electrostatic potential
\begin{align*}
\extpot[\nuc](\el)&:=-\sum_{\nucindex=1}^\nnuc\frac{\charge_\nucindex}{|\el-\nuc_\nucindex|}
\end{align*}
is an external potential, generated by the nuclei, which represents the Coulombic nucleus-electron interactions. The Hartree potential
\begin{align*}
\potH[\dens]:= |\cdot |^{-1}\star\dens
\end{align*}
corresponds to the Coulombic electron-electron interactions. The remaining term, the exchange-correlation potential $\potX+\potC$, is not explicitly known: in the Local-Density Approximation (LDA) introduced by Kohn \&{} Sham in \cite{Kohn1965Self-ConsistentEffects}, for the exchange potential $\potX$ an approximation based on the homogeneous electron gas approximation is chosen \cite{Parr1989Density-FunctionalMolecules}. In this paper, we study a generalisation of this exchange potential, of the form
\begin{align*}
\potX[\dens]:=\lambda\dens^{q-1},
\end{align*}
where $\lambda\in\bb{R}$, $q>1$. Hereafter, we set the so-called correlation potential to zero, namely
\begin{align*}
\potC\equiv 0,
\end{align*}
and write accordingly $\pot_{\text{Hxc}}=\pot_{\text{Hx}}$. In most cases, there is no closed form for the correlation potential, and one has to resort to numerical presentations, which are too complex to investigate in the same manner we handle the other terms. See e.g. \cite{Jerome2015TimeSolutions,Anantharaman2009ExistenceChemistry} and references therein, 
where the case $\potC\not\equiv 0$ is considered in time-independent, resp. specific time-dependent settings.

In the coupling of \eqref{eq:KS} with the equations \eqref{eq:KSNb} describing the nuclear dynamics, which using \eqref{eq:accel} can be written as
\begin{align}
\label{eq:N}
      \ddot{\nuc}&= \acc[\dens](\nuc),
\end{align}
we apply the so-called mean-field, or Ehrenfest dynamics approach, see e.g. \cite{Tully1998MixedDynamics}, \cite[\S 2.3]{Micha2007QuantumSystems}, \cite[\S V]{Agostini2013MixedProcesses} and \cite[\S 2.1]{Crespo-Otero2018RecentDynamics},  based on factorising the total wave function into a product of fast (electronic) and slow (nuclear) particle parts. In this nonadiabatic mixed quantum-classical dynamics method, we use a point-nuclei rather than the Born--Oppenheimer approximation, which would assume some requirements for the system under consideration. This way, we can neglect the quantum nature of the nuclei, since these are much heavier than electrons, and consider them as classical point particles. This mean-field description can be understood as a semi-classical limit of the time-dependent self-consistent field (or Hartree) method, from which the Hamilton--Jacobi equation (equivalent to Newton's law of motion) for the nuclei can be derived. According to this description the nuclei move subject to a single effective potential of Hellman--Feynman type, corresponding to an average over quantum states:
\begin{align*}
\mass_{\nucindex}\acc_{\nucindex}[\dens] (\nuc )=-\nabla_{ \nuc_{\nucindex}}W[\dens](\nuc)\qquad\text{for all $\nucindex$},
\end{align*}
where  
\begin{align}
W[\dens](\nuc)&:= \langle \extpot[\nuc],\dens \rangle_{L^2 (\bb{R}^3)}+W_{\text{nn}}(\nuc),\nonumber\\
\qquad W_{\text{nn}}(\nuc)&:=\tfrac{1}{2}\sum_{\nucindex,\nucindextwo=1,\nucindex\neq \nucindextwo}^{\nnuc}\frac{\charge_{\nucindex} \charge_{\nucindextwo}}{ | \nuc_{\nucindex}-\nuc_{\nucindextwo} |},\label{eq:W2}
\end{align}
describe the interaction of the electrons with the external potential, and the Coulombic internal nuclear interactions by $W_{\text{nn}}$. Note that the exchange term does not appear in the coupling of \eqref{eq:KS} with \eqref{eq:KSNb}, as it does not describe electrostatic interaction, but interactions between the electrons. Also, we note that our equations \eqref{eq:KSN} can be regarded as a Hamiltonian system. The total energy $\energy$ associated with this system is given by
\begin{align*}
\energy[\nuc,{}\wf{}]:=\energy_{\text{kin}} [\nuc,{}\wf{}]+W\big[|\wf|^2\big](\nuc)+\energy_{\text{H}}\big[|\wf|^2\big]+\energy_\text{x}\big[|\wf|^2\big],
\end{align*}
where 
\begin{align*}
\energy_{\text{kin}}[\nuc,{}\wf{}]:=\tfrac{1}{2}\sum_{\nucindex=1}^\nnuc\mass_{\nucindex}\big|\dot{\nuc}_{\nucindex}\big|^2+\tfrac{1}{2}\sum_{\elindex=1}^\nel\int|\nab\wf_\elindex(\el)|^2\rd\el
\end{align*}
is the kinetic energy of the system. The other terms are potential energies:
\begin{align*}
\energy_{\text{H}}\big[|\wf|^2\big]:=\tfrac{1}{2}\iint 
\frac{|\wf(\el)|^2|\wf(\el')|^2}{|\el-\el'|}\rd\el\,\rd\el'
\end{align*}
is the Hartree electrostatic self-repulsion of the KS electron density, and
\begin{align*}
\energy_\text{x}\big[|\wf|^2\big]:=\frac{\lambda}{q}\int|\wf(\el)|^{2q}\rd\el
\end{align*}
is the exchange energy, whose functional derivative coincides with the exchange potential~$\potX$. The total energy $\energy$ as well as $ \|\wf \|_{L^2(\bb{R}^3;\bb{C}^\nel)}$ are quantities which are conserved under the dynamics, as is customary for Hamiltonian systems.

Canc\`es \&{} Le Bris \cite{Cances1999OnDynamics} have considered similar electronic evolution equations coupled with classical nuclear dynamics consistent with the mean-field Ehrenfest approach. They studied a system involving the Hartree--Fock equations:
\begin{subequations}
\label{eq:HFN}
\begin{empheq}[]{alignat=2}
i\dot{\wf}^{\text{HF}}&=\schop^{\text{HF}}[\nuc,\wf^{\text{HF}}]\wf^{\text{HF}},\label{eq:HFNa}\\
\ddot{\nuc}&=\acc[\dens](\nuc),\label{eq:HFNb}\end{empheq}
\end{subequations}
where $\dens=\displaystyle\sum_{\elindex=1}^\nel |\wf^{\text{HF}}_\elindex |^2$, the Hartree--Fock Hamiltonian is defined as
\begin{align*}
\schop^{\text{HF}}[\nuc,\wf^{\text{HF}}]:=-\tfrac{1}{2}\lap+\extpot[\nuc]+\potH[\dens]+\potX^{\text{HF}}[\wf^{\text{HF}}],
\end{align*}
and
\begin{align*}
 (\potX^{\text{HF}}[\wf^{\text{HF}}]\wf^{\text{HF}} )_\elindex&:=-\sum_{\elindextwo=1}^\nel (\overline{\wf^{\text{HF}}_\elindextwo}\wf^{\text{HF}}_\elindex\star |\cdot |^{-1} )\wf^{\text{HF}}_\elindextwo,
\end{align*}
is known as the Hartree--Fock exchange potential. Here, $\wf^{\text{HF}}_\elindex$ are single-particle wave functions. In \cite{Cances1999OnDynamics}, the result of global-in-time existence and uniqueness of solutions to \eqref{eq:HFN} in $H^2$ is based on the celebrated result by Yajima \cite{Yajima1987ExistenceEquations} on the existence of propagators associated with linear, time-dependent Hamiltonians. The proof in \cite{Cances1999OnDynamics} consists of two main steps: a fixed-point argument to show existence of short-time solutions, based on Lipschitz estimates in $H^2 (\bb{R}^3;\bb{C}^\nel )$, and a Gr\"onwall-type argument which relies on energy conservation, conservation of the $L^2 (\bb{R}^3;\bb{C}^\nel )$ norm of $\wf^{\text{HF}}$, and estimates of the solutions $\wf^{\text{HF}}$ in the $H^2 (\bb{R}^3;\bb{C}^\nel )$ norm. 

To the best of our knowledge, since the paper by Canc\`es \&{} Le Bris \cite{Cances1999OnDynamics}, only a few contributions deal with the coupling of a system describing electronic evolution with nuclear dynamics; this is the case, for instance, of \cite{Baudouin2005ExistenceDynamics}, where existence and regularity questions have been studied for a similar system, in the case $\lambda=0$. Considerable attention has also been devoted to Schr\"odinger--Poisson-type equations, which include the Hartree--Fock and the TDKS equations: see for instance \cite{Mauser2001TheEquation,Catto2013ExistencePrinciple,Bokanowski2003OnSystem,Zagatti1992TheEquations,Chadam1975GlobalEquations,Castella1997L2Effects,Anantharaman2009ExistenceChemistry,Jerome2015TimeSolutions,DaPrato1,DaPrato2}. We also mention \cite{Sprengel2019AEquations}, where existence, uniqueness, and regularity questions are investigated for TDKS equations set on bounded space domains, in relation to control problems. None of the contributions listed above have considered the combined nuclear and electronic dynamics as described in our system.
\subsection{Paper outline}
The paper is organised as follows.

In \S\ref{sec:preplemmas}, we recall the relevant results from Yajima \cite{Yajima1987ExistenceEquations} on the construction and properties of a family of propagators
\begin{align*}
U(t,s):L^2 (\bb{R}^3;\bb{C}^\nel )\longrightarrow L^2 (\bb{R}^3;\bb{C}^\nel ),
\end{align*}
with $t,s\in[0,\Theta]$, associated with the linear parts of the KS Hamiltonians $\schop [\nuc(t),\dens ]$ for $t\in[0,\Theta]$, with~$0<\Theta<\infty$, and some results from Canc\`es \&{} Le Bris \cite{Cances1999OnDynamics} on the bounds on the operator norms of these propagators.\\\\

In \S{}\ref{sec:Bel} and \S{}\ref{sec:bnuc} we define bounded regions $\mathcal{B}_{\mathrm{el}} (\tau
 )$ and $\mathcal{B}_{\mathrm{nuc}} (\tau )$, designed to seek solutions to (\ref{eq:KS}, resp. \ref{eq:N}) on a time interval $[0,\tau]$, and the mappings
\begin{align*}
\mathcal{N}:\mathcal{B}_{\mathrm{el}} (\tau )\longrightarrow \mathcal{B}_{\mathrm{nuc}} (\tau )\cap C^2 ([0,\tau];\bb{R}^{3\nnuc} ),\qquad
\mathcal{E}:\mathcal{B}_{\mathrm{nuc}} (\tau )\longrightarrow\mathcal{B}_{\mathrm{el}} (\tau ),
\end{align*}
which connect these solutions.

In \S{}\ref{sec:lip}, in view of a Duhamel-type argument developed in later sections, we state and prove some Lipschitz estimates on the nonlinear mapping 
\begin{align*}
\wf\longmapsto \pot_{\text{Hx}} [ |\wf |^2 ]\wf.
\end{align*}
The restriction $q\geq 7/2$ arises from these estimates.

Next, we prove in \S{}\ref{sec:locex} that for some $\tau>0$ and any fixed $\wf\in\Bel$, the Cauchy problem \eqref{eq:N} has a unique solution $\nuc\in\Bnuc\cap C^2 ([0,\tau];B_\delta (\nuc^0 ) )$, with $B_\delta (\nuc^0 )$ denoting a closed ball of radius $\delta$ centred around $\nuc^0$, and the mapping~$\mathcal{N}[\wf]=\nuc$ is bounded is bounded with respect to the $C^1 ([0,\tau];\bb{R}^{3\nnuc} )$ topology, and continuous as a map from $C^0\big([0,\tau];\ltwoN\big)$ to $C^0\big([0,\tau];\bb{R}^{3\nnuc}\big)$. We construct these solutions as fixed points of the mapping
\begin{align*}
\mathcal{T} [\nuc ](t)=\nuc^0+\vel^0t+\int_0^t(t-\sigma)\acc(\sigma,\nuc(\sigma))\rd \sigma.
\end{align*}
We stress that here $\acc$ depends on $\wf$.

Further, we prove in \S{}\ref{sec:el} that for $q\geq 7/2$, some $\tau>0$ and any fixed~$\nuc\in\Bnuc$, the Cauchy problem \eqref{eq:KS} has a unique solution $\wf\in\Bel$, and the mapping  $\mathcal{E}[\nuc]=\wf$ is bounded and continuous as a map from $C^0\big([0,\tau];\bb{R}^{3\nnuc}\big)$ to $C^0\big([0,\tau];\ltwoN\big)$.. Similarly, solutions are constructed as fixed points of the mapping
\begin{align*}
\mathcal{F}[\wf](t)=U(t,0)\wf^0-i\int_0^tU(t,\sigma)\pot_{\text{Hx}}[\dens]\wf(\sigma)\rd\sigma.
\end{align*}
Using results from \S{}\ref{sec:preplemmas}, \S{}\ref{sec:lip} and Yajima \cite{Yajima1987ExistenceEquations}, we show that fixed points of this mapping are strong solutions to \eqref{eq:KS}.\\\\
We then prove in \S{}\ref{sec:locun} that for $q\geq 7/2$ and some $\tau>0$, the initial-value problem associated with the problem \eqref{eq:KSN} has a solution $ (\wf,\nuc )$ in $\mathcal{X}(\tau)$. To this end, we construct the mapping 
\begin{align*}
\mathcal{K}:\Bnuc\longrightarrow \Bnuc,\qquad \mathcal{K}=\mathcal{I}\circ \mathcal{N}\circ\mathcal{E},
\end{align*}
where 
\begin{align*}
\mathcal{I}:\Bnuc \cap C^2 ([0,\tau];\bb{R}^{3\nnuc} )&\varlonghookrightarrow\Bnuc
\end{align*}
is the inclusion into $\Bnuc;$ we then apply a Schauder-type argument to $\mathcal{K}$, in the spirit of \cite{Cances1999OnDynamics}.  Unlike in \cite{Cances1999OnDynamics}, we equip $\Bnuc$ with a weaker $C^0$-topology, which takes into account nuclear repulsion. The remainder of this section is devoted to uniqueness.

Finally, the Appendix \ref{sec:notation} is devoted to the notation we systematically use, comprising that for the norms on different function spaces, such as $H^2 (\bb{R}^3;\bb{C}^\nel )$ and Lorentz spaces.
\subsection{Related questions}
Theorem \ref{thm:shorttimeexistence} can be generalised to LDA-type nonlinearities which are either sufficiently smooth at the origin $\dens=0$, or enjoy $H^2$-Lipschitz estimates like those obtained in this paper. This is the case, for instance, of $\lambda_1 \dens^{q_1-1}-\lambda_2 \dens^{q_2-1}$ with $q_1,q_2\geq7/2$ and $\lambda_1,\lambda_2>0,$ which share a similar structure with nonlinearities involved in various well-known models in quantum mechanics, such as the Thomas--Fermi--Dirac--Von Weizs\"acker model \cite{Lieb1983DensitySystems}. For this particular example, working with the same functional setting, it would be interesting to explore, for certain ranges of exponents, the occurrence of either a blow-up at finite time in the norm of the solutions or the existence of maximal solutions defined for all $t\geq0$: see \cite{Cazenave1998AnEquations,Cazenave2003SemilinearEquations}.

Also, it would be interesting to identify a functional setting (and a possibly different proof) --- the most natural one would certainly be $H^1$ --- in order to capture the physically relevant exponent $q=4/3,$ which is not covered in the present work. We wonder if a suitable regularisation `at the origin' of the LDA term for $q=4/3$ would allow to cover this case as a result of a limit process. 
\section{Preliminaries}
\label{sec:preplemmas}
The first observation in this section is that the Newton potential
\begin{align}
G [\phi_1,\phi_2 ]:= (\overline{\phi_1}\,\phi_2 )\star |\cdot |^{-1}, \label{eq:mappingG}
\end{align}
solution to 
\begin{align}
-\lap G=4\pi\overline{\phi_1}\phi_2,\label{eq:deltaG}
\end{align}
defines a mapping $H^2\times H^2\longrightarrow W^{2,\infty}.$

\begin{lemma}
\label{lem:G}
For all $i,j\in\{1,2,3\}$ and every $\el\in\bb{R}^3$ it holds that
\begin{alignat}{3}
 |G [\phi_1,\phi_2 ](\el) |&\lesssim{} \rVert\phi_1\rVert_{L^2} \rVert \nab\phi_2 \rVert_{L^2 },\label{eq:Gestimate1}\\
 |\partial_i G [\phi_1,\phi_2 ](\el) |&\lesssim{}
 \rVert\nab\phi_1 \rVert_{L^2 } \rVert \nab\phi_2 \rVert_{L^2 },\label{eq:Gestimate2}\\
 |\partial_{ij} G [\phi_1,\phi_2 ](\el) |&\lesssim{} \rVert\phi_1 \rVert_{H^2} \rVert\phi_2 \rVert_{H^2}.\label{eq:Gestimate3}
\end{alignat}
\end{lemma}
\begin{proof}
By Hardy's inequality \eqref{eq:hardy} and the properties 
\begin{align*}
\partial_iG [\phi_1,\phi_2 ]&= (\overline{\phi_1}
    \phi_2 )\star  (\el_i |\el |^{-3} ),\\
\partial_{ij}G [\phi_1,\phi_2 ]&= [ (\partial_{ i}\overline{\phi_1} )\phi_2+\overline{\phi_1} (\partial_{j}\phi_2 ) ]\star  (\el_i |\el |^{-3} ),
\end{align*}
for all $i,j$ and $\el\in\bb{R}^3$, it holds
\begin{align*}
 |G [\phi_1,\phi_2 ](\el) |&= | \langle\phi_1, |\cdot-\el |^{-1}\phi_2 \rangle_{L^2} |\lesssim{}  \rVert\phi_1 \rVert_{L^2} \rVert \nab\phi_2 \rVert_{L^2 },\\
 |\partial_i G [\phi_1,\phi_2 ](\el) |&\leq  \langle 
 |\cdot-\el |^{-1} |\phi_1 |, |\cdot-\el |^{-1} |\phi_2 | \rangle_{L^2}\lesssim{}
 \rVert\nab\phi_1 \rVert_{L^2 } \rVert \nab\phi_2 \rVert_{L^2 },\\
 |\partial_{ij} G [\phi_1,\phi_2 ](\el) |&\leq
 \langle  |\cdot-\el |^{-1} |\partial_i\phi_1 |, |\cdot-\el |^{-1} |\phi_2 | \rangle_{L^2}+ \langle |\cdot-\el |^{-1} |\phi_1 |, |\cdot-\el |^{-1} |\partial_j\phi_2 | \rangle_{L^2}\nonumber\\
&\lesssim{} \rVert\nab\partial_i\phi_1 \rVert_{L^2 } \rVert\nab\phi_2 \rVert_{L^2 }+ \rVert \nab\phi_1 \rVert_{L^2 } \rVert\nab\partial_j\phi_2 \rVert_{L^2 }\nonumber\\
&\lesssim{} \rVert\phi_1 \rVert_{H^2} \rVert\phi_2 \rVert_{H^2}.
\end{align*}
This concludes the proof.
\end{proof}
The following lemma generalises \cite[Lemma 3]{Cances1999OnDynamics}, and provides us with useful bounds on the  functions  $f_{\nucindex}^{\elindex\elindextwo}:\bb{R}^{3\nnuc}\longrightarrow\bb{C}^3$ defined as
\begin{align*}
f_\nucindex^{\elindex\elindextwo}:=\nabla_{\nuc_{\nucindex}}(\wf_\elindex,\extpot[\nuc]\wf_\elindextwo)_{L^2},
\end{align*}
namely,
\begin{align*}
f_\nucindex^{\elindex\elindextwo}(\nuc)&=-\charge_{\nucindex}\Big(\wf_\elindex,\frac{\cdot-\nuc_{\nucindex}}{|\cdot-\nuc_{\nucindex}|^{3}}\wf_\elindextwo\Big)_{L^2}.
\end{align*}
Note that $f_{\nucindex}^{\elindex\elindextwo}$ effectively only depends on the position $ \nuc_{\nucindex}$ of the $K$-th nucleus, and that
\begin{align*}
\acc^1_\nucindex=-\frac{1}{\mass_{\nucindex}}\sum_{\elindex=1}^\nel f^{\elindex\elindex}_{\nucindex}.
\end{align*}
\begin{lemma}
\label{lem:forces}
For all $\wf_\elindex,\wf_\elindextwo\in H^2$, it holds that 
\begin{align*}
\big\|f_\nucindex^{\elindex\elindextwo}\big\|_{L^{\infty}(\bb{R}^{3\nnuc};\bb{C}^3)}&\lesssim{}
    \|\nab \wf_\elindex\|_{L^2}\|\nab \wf_\elindextwo\|_{L^2},
\end{align*}
and 
\begin{align}
\big\|Df_\nucindex^{\elindex\elindextwo}\big\|_{L^{\infty}(\bb{R}^{3\nnuc};\bb{C}^{3\times 3})}&\lesssim{}\|\wf_\elindex\|_{H^2}\|\wf_\elindextwo\|_{H^2}.\label{eq:forcebounds2}
\end{align}
Here, $D$ is the gradient in $\mathbb R^{3\nnuc}$. In addition, we have that $f_\nucindex^{\elindex\elindextwo}\in W^{1,\infty}\cap C^1$ for all~$\nucindex$.
\end{lemma}
\begin{proof}
By Lemma \ref{lem:G}, $G[\phi_1,\phi_2]\in W^{2,\infty}$ for all $\phi_1,\phi_2\in H^2$. Using 
\begin{align*}
f^{\elindex\elindextwo}_\nucindex(\nuc)=-\charge_{\nucindex}\nab G[\wf_\elindex,\wf_\elindextwo](\nuc_{\nucindex}),
\end{align*}
we get
\begin{align*}
\big\|f^{\elindex\elindextwo}_\nucindex\big\|_{L^{\infty}(\bb{R}^{3\nnuc};\bb{C}^3)}&\lesssim{}\rVert\nab\wf_\elindex\rVert_{L^2}\rVert\nab\wf_\elindextwo\rVert_{L^2},\\
\big\|Df^{\elindex\elindextwo}_\nucindex\big\|_{L^{\infty}(\bb{R}^{3\nnuc};\bb{C}^{3\times 3})}&\lesssim\max_{\nuc_{\nucindex}\in\bb{R}^3}\big\| D^2G[\wf_\elindex,\wf_\elindextwo](\nuc_{\nucindex})\big\|_{\bb{C}^{3\times3}} \lesssim{}\rVert\wf_\elindex\rVert_{H^2}\rVert \wf_\elindextwo\rVert_{H^2}.
\end{align*}
This shows that $f^{\elindex\elindextwo}_\nucindex\in W^{1,\infty}$. By Sobolev's embedding in H\"older spaces, $\overline{\wf_\elindex}\wf_\elindextwo\in C^{0,\alpha}_{\mrm{loc}}$. Using \eqref{eq:mappingG} from Lemma \ref{lem:G} and standard elliptic regularity, it holds that $G[\wf_\elindex,\wf_\elindextwo]\in C^2$, by which~$f^{\elindex\elindextwo}_\nucindex\in C^1$.
\end{proof}
In what follows we recall some results on the existence of the propagator for the linear parts of the Kohn--Sham-type Hamiltonian $\schop [\nuc(t),\dens ]$ for $t\in[0,\Theta]$, with $0<\Theta<\infty$, for a given nuclear configuration $\nuc\in C^1 ([0,\Theta];\bb{R}^{3\nnuc} )$.\\\\
For some $\nuc\in C^1 ([0,\Theta];\bb{R}^{3\nnuc} )$ and $0<\Theta<\infty$ fixed, we consider the family of linear time-dependent Hamiltonians $ \{\schoplin(t),t\in[0,\Theta] \}\subset \mathcal{L} ( \htwoN;\ltwoN ):$
\begin{align}
\label{eq:linHam}
\schoplin(t):=-\tfrac{1}{2}\lap+ \pot(t),
\end{align}
where
\begin{align}
\pot(t,\cdot)&:=\extpot[\nuc(t)],\label{eq:calV}
\end{align}
Note that $\schoplin(t)$ is the linear part of $\schop [\nuc(t),\dens ]$, and that for any fixed $t$ is a self-adjoint operator on $\ltwoN$. We emphasise that these expressions depend on the time evolution of the nuclear configuration $\nuc$. This family of Hamiltonians is naturally associated with the Cauchy problem
\begin{align*}
i\dot{\wf}= \schoplin(t)\wf,\qquad \wf(s)=\wf^0,
\end{align*}
on a time interval $[0,\Theta]$, for some $s\in[0,\Theta]$. Equivalently, we can formulate the above as an integral equation
\begin{align}
\label{eq:IE}
\wf(t)=U_0(t-s)\wf^0-i \int_s^tU_0(t-\sigma)\pot(\sigma)\wf(\sigma)\rd \sigma,
\end{align}
where
\begin{align*}
U_0(t):=\exp (it\lap/2)
\end{align*}
is the free propagator (i.e., the propagator for the free particle), which is an evolution operator on $ \htwoN$. The following lemma is in the spirit of \cite[Lemma 4]{Cances1999OnDynamics}, which in turn is based on \cite[Cor. 1.2. (1)--(2)--(4), Thm. 1.1. (2) \& Thm. 1.3. (5)--(6)]{Yajima1987ExistenceEquations}.
\begin{lemma}
\label{lem:evolutionoperators}
For the family of Hamiltonians $ \{\schoplin(t),t\in[0,\Theta] \}$, there exists a unique family of linear evolution operators 
\begin{align*}
U(t,s):\ltwoN\longrightarrow\ltwoN,\qquad t,s\in[0,\Theta],
\end{align*}
such that 
\begin{align*}
    \wf(t):=U(t,s)\wf^0
\end{align*}
solves \eqref{eq:IE} on $[0,\Theta]$ for all $\wf^0\in \htwoN$, with 
\begin{align*}
 \|\wf(t) \|_{\ltwoN}= \|\wf^0 \|_{\ltwoN}
\end{align*}
for all $t\in[0,\Theta]$. Moreover, this family enjoys the following properties: 
\begin{enumerate}[(i)]
\item $U(t,s)U(s,r)=U(t,r)$ for all $t,s,r\in[0,\Theta]$.
\item $U(t,t)=\mrm{Id}$ for all $t\in[0,\Theta]$.
\item $U(t,s)$ is a unitary operator on $\ltwoN$ for all $t,s\in[0,\Theta]$:
\begin{align*}
 \|U(t,s)\wf \|_{\ltwoN}= \|\wf \|_{\ltwoN}.
\end{align*}
\item For all $f\in\ltwoN$, 
$ ((t,s)\longmapsto U(t,s)f ): [0,\Theta]^2\longrightarrow \ltwoN$ is a continuous mapping.
\item $U(t,s)\in \mathcal{L} ( \htwoN )$ for all $(t,s)\in[0,\Theta]^2$.
\item For all $f\in \htwoN$, 
$ ((t,s)\longmapsto U(t,s)f ): [0,\Theta]^2\longrightarrow  \htwoN$ is a continuous mapping.
\item For all $f\in \htwoN$, the mapping $(t,s)\longmapsto U(t,s)f$ is an element in $C^1 ([0,\Theta]^2;\ltwoN )$, and the following equations hold in $\ltwoN$:
\begin{align*}
i\frac{\partial}{\partial t} (U(t,s)f )&=\schoplin(t)U(t,s)f,\\
i\frac{\partial}{\partial s} (U(t,s)f )&=-U(t,s)\schoplin(s)f.
\end{align*}
\item For all $\gamma>0$, there is a constant $B_{\Theta,\gamma}$ of the form
\begin{align*}
B_{\Theta,\gamma}=A_\gamma^{1+C_\gamma\Theta},\qquad A_\gamma,C_\gamma>0,
\end{align*}
such that if
\begin{align*}
 \|\dot{\nuc} \|_{C^0 ([0,\Theta];\bb{R}^{3\nnuc} )}\leq \gamma,
\end{align*}
then for all $t,s\in[0,\Theta]$
\begin{align*}
 \|U(t,s) \|_{\mathcal{L} (\htwoN)}\leq B_{\Theta,\gamma}.
\end{align*}
\end{enumerate}
\end{lemma}
\begin{proof}
The result in the case $\nel=\nnuc=1$ has been proved in \cite[Lemma 4]{Cances1999OnDynamics}. We observe that the argument in \cite{Cances1999OnDynamics} is robust enough to be easily adapted to our more general context of arbitrary $\nel,\nnuc\in\bb{N}$. Indeed, since the linear Hamiltonians $\schoplin(t)$ do not depend on an electronic configuration~$\wf$, and act on every element $\wf_\elindex$ independently, the result for general $\nel$ follows from the case $\nel=1$. In particular, properties (i)---(vii) can be justified with an obvious adaptation of the case $\nel=1$ proved in \cite[Cor. 1.2. (1)--(2)--(4), Thm. 1.1. (2) \& Thm. 1.3. (5)--(6)]{Yajima1987ExistenceEquations}. We note that property (viii) can be also justified arguing exactly as for the case $\nel=\nnuc=1$ in \cite[Lemma 4]{Cances1999OnDynamics}, observing that our additional terms in the expression of $\pot$ can be estimated in the same way.
\end{proof}
\section{Definition of the electronic feasible region \texorpdfstring{$\mathcal{B}_{\mathrm{el}}$}{Bel}}
\label{sec:Bel}
Let $\tau>0$ be finite and define
\begin{align}
\gamma:= |\vel^0 |+1,\label{eq:gamma}
\end{align}
where the term ``$+1$'' allows us to cover the case $\vel^0=0$.
Let us consider $B_{\tau,\gamma}$ as given in Lemma \ref{lem:evolutionoperators} with $\Theta=\tau,$ and where $\gamma$ is as above. We can therefore define the radius
\begin{align*}
\alpha(\tau)&:=2B_{\tau,\gamma} \|\wf^0 \|_{ \htwoN}
\end{align*}
for the ball centred around the initial configuration $\wf^0\in \htwoN$:
\begin{align*}
B_\alpha (\wf^0 )&= \{\wf\in \htwoN \rvert \|\wf-\wf^0 \|_{ \htwoN}\leq \alpha \}.
\end{align*}
Finally, let us define the electronic feasible region for the time interval $[0,\tau]$ as 
\begin{align*}
\mathcal{B}_{\mathrm{el}} (\tau
 )&:= \{\wf \in C^1 ([0,\tau];\ltwoN )\cap C^0 ([0,\tau];B_{\alpha} (\wf^0 ) )\,\, \rvert
\,\,\wf(0)=\wf^0 \},
\end{align*}
equipped with the $C^0 ([0,\tau];\ltwoN )$ norm, which is designed to contain solutions $\wf$ to the Cauchy problem associated with \eqref{eq:KS} with $\wf(0)=\wf^0$ on the time interval $[0,\tau]$, which we may call feasible electronic configurations.
\section{Definition of the nuclear feasible region \texorpdfstring{$\mathcal{B}_{\mathrm{nuc}}$}{Bnuc}}
\label{sec:bnuc}
For all $0<\varepsilon<\min_{\nucindex\neq \nucindextwo} \{ | \nuc_{\nucindex}^0-\nuc_{\nucindextwo}^0 | \}$, we set \begin{align*}\delta(\tau):=\frac{\min_{\nucindex\neq \nucindextwo} \{ | \nuc_{\nucindex}^0-\nuc_{\nucindextwo}^0 | \}-\min \{ \delta_{\textrm{rep}}(\tau),\varepsilon \}}{2}>0,\end{align*}
where 
\begin{align*}
\delta_{\textrm{rep}}(\tau)&:=\biggr[\biggr(\sum_{\nucindex=1}^\nnuc\mass_{\nucindex}\big|\vel^0_\nucindex\big|^2+
\sum_{\subalign{\nucindex,&\nucindextwo=1,\\\nucindextwo&\neq\nucindex}}^\nnuc\frac{\charge_{\nucindex} \charge_{\nucindextwo}}{\big|\nuc_{\nucindex}^0-\nuc_\nucindextwo^0\big|}\biggr)\mrm{e}^\tau\nonumber\\&\qquad +
16\sum_{\nucindex=1}^\nnuc\frac{\charge_\nucindex^2}{\mass_\nucindex}\|{}\wf{}\|_{C^0([0,\tau];\honeN)}^2\big(\mrm{e}^\tau-1\big)\biggr]^{-1}
\end{align*}
arises from a repulsion argument, given in the lemma below. Note that $\varepsilon>0$ ensures the strict positivity of $\delta(\tau)$, which defines the radius for the ball centred around the initial configuration $\nuc^0\in\bb{R}^{3\nnuc}$, with $\nuc^0_{\nucindex}\neq \nuc^0_{\nucindextwo}$ for $\nucindex\neq \nucindextwo$:
\begin{align*}
B_\delta(\nuc^0)&= \{\nuc\in\bb{R}^{3\nnuc} \rvert |\nuc-\nuc^0 |\leq \delta \}.
\end{align*}
Then, by the triangle inequality, for all $\nuc\in B_\delta (\nuc^0 )$ and $\nucindex\neq\nucindextwo$,  it holds that
\begin{align*}
&|\nuc_{\nucindex}-\nuc_\nucindextwo|
\geq
\min_{\nucindex'\neq \nucindextwo'}\big|\nuc_{\nucindex'}^0-\nuc_{\nucindextwo'}^0\big|
-2\big|\nuc-\nuc^0\big|\nonumber\\
&\qquad \geq \min_{\nucindex'\neq \nucindextwo'}\big|\nuc_{\nucindex'}^0-\nuc_{\nucindextwo'}^0\big|
-2\delta(\tau)=\min\{ \delta_{\textrm{rep}}(\tau),\varepsilon\}>0.
\end{align*}
We define the nuclear feasible region for the time interval $[0,\tau]$ as
\begin{align*}
\mathcal{B}_{\mathrm{nuc}} (\tau )
&:= \{\nuc\in C^1 ([0,\tau];B_{\delta} (\nuc^0 ) )
\, \rvert\,\nuc(0)=\nuc^0,\dot{\nuc}(0)=\vel^0, \|\dot{\nuc} \|_{C^0 ([0,\tau];\bb{R}^{3\nnuc} )}\leq \gamma \}
\end{align*}
with $\gamma$ as in \eqref{eq:gamma}. This region is equipped with the $C^0 ([0,\tau];\mathbb{R}^{3\nnuc} )$
topology, and is designed to contain short-time solutions $\nuc$ to the Cauchy problem associated with \eqref{eq:N} with $\nuc(0)=\nuc^0,\dot{\nuc}(0)=\vel^0$ on the interval $[0,\tau]$, which we call feasible nuclear configurations.

This definition of $\delta(\tau)$ is suggested by an a priori lower bound on the nuclear distances $ | \nuc_{\nucindex}(t)-\nuc_{\nucindextwo}(t) |$, $\nucindex\neq \nucindextwo$, which is based on Gr\"onwall's lemma. In fact, we have the following

\begin{lemma}\label{lem:deltaGW2}
Fix ${}\wf{}\in C^0\big([0,\tau];\honeN\big)$, and $\nuc^0\in\bb{R}^{3\nnuc}$ such that $\nuc_\nucindex^0\neq\nuc_\nucindextwo^0$ for $\nucindex\neq\nucindextwo$. Let $\nuc$ solve \eqref{eq:N}, and $\nuc(0)=\nuc^0$. Then, for all $t\in[0,\tau]$ and $\nucindex\neq\nucindextwo$
\begin{align*}
&|\nuc_{\nucindex}(t)-\nuc_\nucindextwo(t)|\geq\delta_{\textrm{rep}}(\tau).
\end{align*}

\end{lemma}
\begin{proof}
Writing the momenta $\mom_\nucindex:=\mass_{\nucindex}\dot{\nuc}_\nucindex$, we define the classical reduced Hamiltonian
\begin{align*}
\mathcal{H}_{\text{nn}}(\nuc,\mom):=\tfrac{1}{2}\sum_{\nucindex=1}^\nnuc\frac{|\mom_\nucindex|^2}{\mass_{\nucindex}}+W_{\text{nn}}(\nuc)
\end{align*}
with $W_{\text{nn}}$ as in \eqref{eq:W2}. Fix ${}\wf{}\in C^0\big([0,\tau];\honeN\big)$. Now,
\begin{align}
&\frac{\mrm{d}}{\mrm{d} t}[\mathcal{H}_{\text{nn}}(\nuc,\mom)]=
\sum_{\nucindex=1}^\nnuc\big[\nabla_{\nuc_{\nucindex}}\mathcal{H}_{\text{nn}}(\nuc,\mom)\cdot\dot{\nuc}_{\nucindex}
+\nabla_{\mom_\nucindex}\mathcal{H}_{\text{nn}}(\nuc,\mom)\cdot\dot{\mom}_\nucindex
\big]\nonumber\\
&=
\sum_{\nucindex=1}^\nnuc\frac{\mom_\nucindex}{\mass_{\nucindex}}\cdot\{\nabla_{\nuc_{\nucindex}}[W_{\text{nn}}(\nuc)]
+\mass_{\nucindex}\ddot{\nuc}_\nucindex
\}\nonumber\\
&\stackrel{\eqref{eq:N}}{=}-\sum_{\nucindex=1}^\nnuc\frac{\mom_\nucindex}{\mass_{\nucindex}}\cdot
( \nabla_{\nuc_{\nucindex}}\extpot[\nuc],\dens)_{L^2}\nonumber\\
&\leq\sum_{\nucindex=1}^\nnuc\frac{1}{2\mass_{\nucindex}}\big[|\mom_\nucindex|^2
+|( \nabla_{\nuc_{\nucindex}}\extpot[\nuc],\dens)_{L^2}|^2\big]\nonumber\\
&\leq
\mathcal{H}_{\text{nn}}(\nuc,\mom)+
8\sum_{\nucindex=1}^\nnuc\frac{\charge_\nucindex^2}{\mass_\nucindex}\|{}\wf{}\|_{C^0([0,\tau];\honeN)}^2,
\label{eq:gronwall1}\end{align}
by which, using Gr\"onwall's inequality, on $[0,\tau]$
\begin{align}
&\mathcal{H}_{\text{nn}}(\nuc,\mom)\leq \nonumber\\
&\leq \tfrac{1}{2}\mrm{e}^\tau\biggr[\sum_{\nucindex=1}^\nnuc\mass_{\nucindex}\big|\vel^0_\nucindex\big|^2+
\sum_{\subalign{\nucindex,&\nucindextwo=1,\\\nucindextwo&\neq\nucindex}}^\nnuc\frac{\charge_{\nucindex} \charge_{\nucindextwo}}{\big|\nuc_{\nucindex}^0-\nuc_\nucindextwo^0\big|} +
16\sum_{\nucindex=1}^\nnuc\frac{\charge_\nucindex^2}{\mass_\nucindex}\|{}\wf{}\|_{C^0([0,\tau];\honeN)}^2\biggr]\nonumber
\\
&\qquad -8\sum_{\nucindex=1}^\nnuc\frac{\charge_\nucindex^2}{\mass_\nucindex}\|{}\wf{}\|_{C^0([0,\tau];\honeN)}^2.\label{eq:gronwall2h2}
\end{align}
In \eqref{eq:gronwall1}, we use that, by Hardy's inequality, for all $\nucindex=1,\ldots,\nnuc$
\begin{align*}
&|(\nabla_{\nuc_{\nucindex}}\extpot[\nuc],\dens)_{L^2}|=\nonumber\\
&=\Big|\Big( -\charge_{\nucindex}\frac{\cdot-\nuc_{\nucindex}}{|\cdot-\nuc_{\nucindex}|^3},\dens\Big)_{L^2}\Big|\leq \charge_{\nucindex}\sum_{\elindex=1}^\nel\big\||\cdot-\nuc_{\nucindex}|^{-1}\wf_\elindex\big\|_{L^2}\\
&\leq 
2\charge_{\nucindex}\sum_{\elindex=1}^\nel\|\nab\wf_\elindex\|_{L^2}\leq 
2\sqrt{2}\charge_{\nucindex}\|\nab{}\wf{}\|_{\ltwoN}\leq 2\sqrt{2}\charge_{\nucindex}\|{}\wf{}\|_{\honeN}.\nonumber
\end{align*}
Estimate \eqref{eq:gronwall2h2} is enough to conclude, as for all $\nucindex\neq\nucindextwo$ we have
\begin{align*}
\tfrac{1}{2}\frac{1}{|\nuc_{\nucindex}-\nuc_\nucindextwo|}\leq \frac{1}{\charge_{\nucindex}\charge_{\nucindextwo}}W_{\text{nn}}(\nuc)\leq W_{\text{nn}}(\nuc)\leq \mathcal{H}_{\text{nn}}(\nuc,\mom).
\end{align*}

\end{proof}
\begin{remark}
\label{rem:GW}
A similar argument yields an a priori estimate of the nuclear velocity $\dot{\nuc}$.
\end{remark}
\section{Lipschitz estimates}
\label{sec:lip}
In this section, we obtain Lipschitz estimates on the mapping $\wf\longmapsto \pot_{\text{Hx}}[\dens]\wf:= (\potH[\dens]+\potX[\dens] )\wf$.
\begin{lemma}[Lipschitz estimates on the Hartree term]
\label{lem:lipschitzconvolution}
For all $\wf,\wf'\in \honeN$, with~$\dens':= |\wf' |^2$,
\begin{align}
& \|\potH[\dens]\wf-\potH [\dens' ]\wf' \|_{\ltwoN}\lesssim{} \sqrt{\nel} \|\wf-\wf' \|_{\ltwoN}\times\nonumber\\
&\quad\times \biggr[\sum_{\elindex=1}^\nel ( \|\nab\wf_\elindex \|_{L^2 (\bb{R}^3;\bb{C}^3 )}+
 \|\nab\wf'_\elindex \|_{L^2 (\bb{R}^3;\bb{C}^3 )} )
 \|\wf' \|_{\ltwoN}+
\sum_{\elindextwo=1}^\nel \|\wf_\elindextwo \|_{L^2} \|\nab\wf_\elindextwo \|_{L^2 (\bb{R}^3;\bb{C}^3 )} \biggr].\label{eq:A}
\end{align}
Moreover, for all $\wf,\wf'\in  \htwoN$
\begin{align}
 \|\potH[\dens]\wf \|_{\htwoN}
&\lesssim{} \sqrt{\nel}
\sum_{\elindex=1}^\nel \|\wf_\elindex \|_{H^1}^2
 \|\wf \|_{\htwoN},\label{eq:B}\\
 \|\potH[\dens]\wf-\potH [\dens' ]\wf' \|_{\htwoN}
&\lesssim{}\sqrt{\nel} \|\wf-\wf' \|_{\htwoN}\times \nonumber\\
&\qquad\sum_{\elindex=1}^\nel [ ( \|\wf_\elindex \|_{H^1}+
 \|\wf_\elindex' \|_{H^1} )
 \|\wf' \|_{\htwoN}+ \|\wf_\elindex \|^2_{H^1} ].
\label{eq:C}
\end{align}
\end{lemma}
\begin{proof}
\textit{Proof of \eqref{eq:A}.}\\
By adding and subtracting the term $ ( |\wf'_\elindextwo |^2\star |\cdot |^{-1} )\wf_\elindex'$, we can write for all $\elindex=1,\ldots,\nel$
\begin{align}
& \| (\potH[\dens]\wf-\potH [\dens' ]\wf' )_\elindex \|_{L^2}\leq\nonumber\\
&\qquad\sum_{\elindextwo=1}^\nel
 [\underbrace{ \| ( |\wf_\elindextwo |^2\star |\cdot |^{-1} ) (\wf_\elindex-\wf_\elindex' ) \|_{L^2}}_{=:\,\mrm{(I)}}+
\underbrace{ \| ( ( |\wf_\elindextwo |^2- |\wf_\elindextwo' |^2 )\star |\cdot |^{-1} )\wf_\elindex' \|_{L^2}}_{=:\,\mrm{(II)}}
 ].\label{eq:toA}
\end{align}
Using the Cauchy--Schwarz inequality in (\ref{eq:star11},\ref{eq:star31}), Hardy's inequality in (\ref{eq:star21},\ref{eq:star41}), and the triangle inequality in \eqref{eq:star41},
\begin{align}
\mrm{(I)}&\leq  \| |\wf_\elindextwo |^2\star |\cdot |^{-1} \|_{L^{\infty}} \|\wf_\elindex-\wf_\elindex' \|_{L^2}\leq \underset{\el\in\bb{R}^3}{\mrm{esssup}} \{ | \langle  |\wf_\elindextwo |, |\cdot-\el |^{-1} |\wf_\elindextwo | \rangle_{L^2} | \} \|\wf-\wf' \|_{\ltwoN}\nonumber\\
&\leq
\underset{\el\in\bb{R}^3}{\mrm{esssup}} \{ \|
\wf_\elindextwo \|_{L^2} \| |\cdot-\el |^{-1}\wf_\elindextwo \|_{L^2} \} \|\wf-\wf' \|_{\ltwoN}\label{eq:star11}\\
&\lesssim{} \|\wf_\elindextwo \|_{L^2} \|\nab\wf_\elindextwo \|_{L^2 (\bb{R}^3;\bb{C}^3 )} \|\wf-\wf' \|_{\ltwoN},\label{eq:star21}
\end{align}
and
\begin{align}
\mrm{(II)}&\leq \| ( |\wf_\elindextwo |^2- |\wf_\elindextwo' |^2 )\star |\cdot |^{-1} \|_{L^{\infty}} \|\wf_\elindex' \|_{L^2}\nonumber\\
&\leq \underset{\el\in\bb{R}^3}{\mrm{esssup}} \{ | \langle  |\wf_\elindextwo |- |\wf_\elindextwo' |, |\cdot-\el |^{-1} ( |\wf_\elindextwo |+ |\wf_\elindextwo' | ) \rangle_{L^2} | \} \|\wf' \|_{\ltwoN}\nonumber\\
&\leq
\underset{\el\in\bb{R}^3}{\mrm{esssup}} \{
 \|
 |\wf_\elindextwo |- |\wf_\elindextwo' | \|_{L^2}
 ( \| |\cdot-\el |^{-1}\wf_\elindextwo \|_{L^2}+
 \| |\cdot-\el |^{-1}\wf_\elindextwo' \|_{L^2} ) \} \|\wf' \|_{\ltwoN}\label{eq:star31}\\
&\lesssim{}
 ( \|\nab\wf_\elindextwo \|_{L^2 (\bb{R}^3;\bb{C}^3 )}+
 \|\nab\wf_\elindextwo' \|_{L^2 (\bb{R}^3;\bb{C}^3 )}
 ) \|\wf' \|_{\ltwoN} \|\wf_\elindextwo-\wf_\elindextwo' \|_{L^2}\label{eq:star41}\\
&\lesssim{} ( \|\nab\wf_\elindextwo \|_{L^2 (\bb{R}^3;\bb{C}^3 )}+
 \|\nab\wf_\elindextwo' \|_{L^2 (\bb{R}^3;\bb{C}^3 )}
 ) \|\wf' \|_{\ltwoN} \|\wf-\wf' \|_{\ltwoN}.\label{eq:alsotoA}
\end{align}
Combining (\ref{eq:star21},\ref{eq:alsotoA}) with \eqref{eq:toA} we get \eqref{eq:A}.\\
\\
\textit{Proof of \eqref{eq:B}.}\\
With $G$ as in \eqref{eq:mappingG}, we can write 
\begin{align*}
\potH[\dens]\wf=\sum_{\elindex=1}^\nel G [\wf_\elindex,\wf_\elindex ]\wf.
\end{align*}
Note that for all $\phi_1,\phi_2,\phi_3\in H^2$
\begin{align*}
\lap (G [\phi_1,\phi_2 ]\phi_3 )\stackrel{\eqref{eq:deltaG}}{=}G [\phi_1,\phi_2 ]\lap\phi_3+2\nab (G [\phi_1,\phi_2 ] )\cdot\nab\phi_3-4\pi\overline{\phi_1}\phi_2\phi_3.
\end{align*}
Using (\ref{eq:Gestimate1},\ref{eq:Gestimate2}) from Lemma \ref{lem:G}, we obtain
\begin{align*}
 \|G [\phi_1,\phi_2 ]\phi_3 \|_{L^2}&\leq  \|G [\phi_1,\phi_2 ] \|_{L^{\infty}} \|\phi_3 \|_{L^2}\lesssim{} \|\phi_1 \|_{H^1} \|\phi_2 \|_{H^1} \|\phi_3 \|_{H^1},\\
 \|G [\phi_1,\phi_2 ]\lap\phi_3 \|_{L^2}&\leq \|G [\phi_1,\phi_2 ] \|_{L^{\infty}} \|\lap\phi_3 \|_{L^2}\lesssim{}  \|\phi_1 \|_{H^1} \|\phi_2 \|_{H^1} \|\phi_3 \|_{H^2},\\
 \|\nab (G [\phi_1,\phi_2 ] )\cdot\nab\phi_3 \|_{L^2}&\leq  \|G [\phi_1,\phi_2 ] \|_{W^{1,\infty}} \|\nab\phi_3 \|_{L^2 (\bb{R}^3;\bb{C}^3 )}\lesssim{}  \|\phi_1 \|_{H^1} \|\phi_2 \|_{H^1} \|\phi_3 \|_{H^1}.
\end{align*}
On the other hand, by H\"older's and Sobolev's inequalities,
\begin{align*}
 \|\overline{\phi_1}\phi_2\phi_3 \|_{L^2}&\leq  \|\phi_1 \|_{L^6} \|\phi_2 \|_{L^6} \|\phi_3 \|_{L^6}\lesssim{}  \|\phi_1 \|_{H^1} \|\phi_2 \|_{H^1} \|\phi_3 \|_{H^1}.
\end{align*}
This gives for all $\elindex=1,\ldots,\nel$
\begin{align*}
 \| (\potH[\dens]\wf )_\elindex \|_{H^2}\leq 
\sum_{\elindextwo=1}^\nel \|G [\wf_\elindextwo,\wf_\elindextwo ]\wf_\elindex \|_{H^2}\lesssim{}\sum_{\elindextwo=1}^\nel \|\wf_\elindextwo \|_{H^1}^2 \|\wf \|_{\htwoN}.
\end{align*}
Combining these estimates \eqref{eq:B} follows.\\\\
\textit{Proof of \eqref{eq:C}.}\\
As in the proof of \eqref{eq:A}, we bound for all $\elindex=1,\ldots,\nel$
\begin{align*}
& \|\lap(\potH[\dens]\wf-\potH [\dens' ]\wf' )_\elindex \|_{L^2}\leq\\
&\qquad\sum_{\elindextwo=1}^\nel \{\underbrace{ \|\lap[G [\wf_\elindextwo,\wf_\elindextwo ] (\wf_\elindex-\wf_\elindex' ) ] \|_{L^2}}_{=:\,\mrm{(I)}}
+\underbrace{ \|\lap[G [ |\wf_\elindextwo |+ |\wf_\elindextwo' |, |\wf_\elindextwo |- |\wf_\elindextwo' | ]\wf_\elindex' ] \|_{L^2}}_{=:\,\mrm{(II)}} 
\}.\nonumber
\end{align*}
As for \eqref{eq:B}, we can bound (I) and (II) using $(\phi_1,\phi_2,\phi_3)=(\wf_\elindextwo,\wf_\elindextwo,\wf_\elindex-\wf_\elindex')$ for (I) and $(\phi_1,\phi_2,\phi_3)= ( |\wf_\elindextwo |+ |\wf_\elindextwo' |, |\wf_\elindextwo |- |\wf_\elindextwo' |,\wf_\elindex )$ for (II). Hence, by the triangle inequality,
\begin{align}
& \|\lap(\potH[\dens]\wf-\potH [\dens' ]\wf' ) \|_{\ltwoN}\lesssim{}\mrm{(B)}:=\nonumber
\\
&\quad = \sqrt{\nel}
\sum_{\elindex=1}^\nel \big[ \|\wf_\elindex \|^2_{H^1}
+ (
 \|\wf_\elindex \|_{H^1}
+
 \|\wf'_\elindex \|_{H^1}
 ) \|\wf' \|_{ \htwoN}
 \big] \|\wf-\wf' \|_{ \htwoN}.
\label{eq:deltal2nb}
\end{align}
On the other hand, by \eqref{eq:A} we also have
\begin{align}
 \|\potH[\dens]\wf-\potH [\dens' ]\wf' \|_{\ltwoN}\lesssim{}\mrm{(B)}.
\label{eq:l2nb}
\end{align}
Hence, by (\ref{eq:deltal2nb}, \ref{eq:l2nb}), \eqref{eq:C} immediately follows.
\end{proof}
\noindent By Cauchy--Schwarz inequality, for all $\bb{C}^\nel$-valued functions $\wf,\wf'$ we have
\begin{align}
\label{eq:ell2bound}
 |\wf\cdot\nab\wf' |
\leq
 |\wf | |\nab\wf' |.
\end{align}
\begin{lemma}[Mean-value estimates for the density]
\label{lem:MVE}
For all $a\geq1/2$, we have
\begin{align}
 |\dens^a-\dens'^a |\lesssim_{a}  ( \|\dens \|^{a-1/2}_{L^{\infty}}+ \|\dens' \|^{a-1/2}_{L^{\infty}} ) |\wf-\wf' |.\label{eq:MVE}
\end{align}
\end{lemma}
\begin{proof}
By the fundamental theorem of calculus
\begin{align*}
 |\dens^a-\dens'^a |&= | |\wf |^{2a}- |\wf' |^{2a} |= \bigg|\int_0^1\frac{\mrm{d}}{\mrm{d}t} [ |\wf'+t (\wf-\wf' ) |^{2a} ]\mrm{d}t \bigg|\\
&\lesssim_{a}  ( |\wf |+ |\wf' | )^{2a-1} |\wf-\wf' |\lesssim_{a} ( |\wf |^{2a-1}+ |\wf' |^{2a-1} ) |\wf-\wf' |\nonumber\\
&=
 (\dens^{a-1/2}+\dens'^{a-1/2} ) |\wf-\wf' |,\nonumber
\end{align*}
which yields \eqref{eq:MVE}.
\end{proof}
\begin{lemma}[Mean-value estimates for the density gradient]
\label{lem:MVEnablarho}
For all $b\geq3/2$, it holds that 
\begin{align}
 |\nab (\dens^{b} )-\nab (\dens'^{b} ) |
\lesssim_{b}\,& (Q_1 |\nab\wf |+Q_2 |\nab\wf' | ) |\wf-\wf' |+Q_3 |\nab\wf-\nab\wf' |,\label{eq:MVEnablarho}
\end{align}
where
\begin{align*}
Q_1=\dens^{b-1},\qquad Q_2= (\dens^{b-3/2}+\dens'^{b-3/2} )\dens'^{1/2},\qquad Q_3=\dens^{b-1}\dens'^{1/2}.
\end{align*}
\end{lemma}
\begin{proof}
Using $\nab\dens=\nab\wf\cdot \overline{\wf}+\wf\cdot\nab (\overline{\wf} )$ and \eqref{eq:ell2bound} for the pair $(\wf,\overline{\wf})$,
\begin{align}
\label{eq:nablarho}
 |\nab\dens |\lesssim{}\dens^{1/2} |\nab\wf |.
\end{align}
Since
\begin{align*}
 |\nab (\dens^{b} ) |\lesssim_{b}\dens^{b-1} |\nab\dens |,
\end{align*}
adding and subtracting the term $\dens^{b-1}\nab\dens'$, we get for all $b\geq 1$
\begin{align*}
 |\nab (\dens^{b} )-\nab (\dens'^{b} ) |\lesssim_{b}\dens^{b-1} |\nab\dens-\nab\dens' |+ |\dens^{b-1}-\dens'^{b-1} | |\nab\dens' |.
\end{align*}
By adding and subtracting $\wf'\cdot\nab (\overline{\wf} )$ and $\overline{\wf'}\cdot\nab\wf$, and using \eqref{eq:ell2bound} for the pairs~$ (\wf-\wf',\overline{\wf} )$, $ (\wf',\overline{\wf-\wf'} )$, $ (\overline{\wf-\wf'},\wf )$ and  $ (\overline{\wf'},\wf-\wf' )$, we get
\begin{align}
\nonumber
 |\nab\dens-\nab\dens' |&=
 |\wf\cdot\nab (\overline{\wf} )-\wf'\cdot\nab (\overline{\wf'} )+\overline{\wf}\cdot\nab\wf-\overline{\wf'}\cdot\nab (\wf' ) |\\
\nonumber
&\leq | (\wf-\wf' )\cdot\nab (\overline{\wf} ) |+
 |\wf'\cdot\nab (\overline{\wf-\wf'} ) |\\
&\qquad \nonumber+ | (\overline{\wf-\wf'} )\cdot\nab \wf  |+ |\overline{\wf'}\cdot\nab (\wf-\wf' ) |\\
&\lesssim{} |\nab\wf | |\wf-\wf' |+\dens'^{1/2} |\nab\wf-\nab\wf' |.
\label{eq:nablarhodiff}
\end{align}
By \eqref{eq:MVE} with $a=b-1\geq 1/2$ and \eqref{eq:nablarho} for $\dens'$, we get
\begin{align*}
 |\dens^{b-1}-\dens'^{b-1} | |\nab\dens' |\lesssim_{b} (\dens^{b-3/2}+\dens'^{b-3/2} )\dens'^{1/2} |\nab\wf' | |\wf-\wf' |.
\end{align*}
Putting these estimates together gives \eqref{eq:MVEnablarho}.
\end{proof}
\begin{lemma}[Lipschitz estimates on the local nonlinearity]
\label{lem:lipschitzlda}
Let $q\in[1,\infty)$, and $\lambda\in\bb{R}$. For any fixed $p\in[1,\infty]$ and for all $\wf,\wf'\in  \htwoN\cap  L^p(\bb{R}^3;\bb{C}^\nel),$ it holds that 
\begin{align}
 \|\potX[\dens]\wf-\potX [\dens' ]\wf' \|_{L^p(\bb{R}^3;\bb{C}^\nel)}
\lesssim_{q,\lambda}
\sum_{\elindex=1}^\nel \Big[ \|\wf_\elindex \|_{H^2}^{2(q-1)}+ \|\wf_\elindex' \|_{H^2}^{2(q-1)} \Big]
 \|\wf-\wf' \|_{ L^p(\bb{R}^3;\bb{C}^\nel)}.\label{eq:D}
\end{align}
Moreover, for all $q\geq 7/2$ and any $\lambda\in\bb{R}$, we have that
\begin{align}
&\|\potX[\dens]\wf-\potX [\dens' ]\wf' \|_{\htwoN}
\leq\nonumber\\
&\quad \leq \mathcal{L}_{q,\lambda} (\max \{ \|\wf \|_{\htwoN}, \|\wf' \|_{\htwoN} \} )
 \|\wf-\wf' \|_{\htwoN},\label{eq:E}
\end{align}
where $\mathcal{L}_{q,\lambda}:\bb{R}_0^+\longrightarrow\bb{R}_0^+$ is a polynomial function which vanishes at the origin.
\end{lemma}
\begin{proof}
\textit{Proof of \eqref{eq:D}.}\\ 
By the fundamental theorem of calculus,
\begin{align*}
 |\potX[\dens]\wf-\potX [\dens' ]\wf' |&=|\lambda| | |\wf |^{2(q-1)}\wf- |\wf' |^{2(q-1)}\wf' |\\
&\lesssim_{\lambda} \Big|\int_0^1\frac{\mrm{d}}{\mrm{d}t} \big[ |\wf'+t (\wf-\wf' ) |^{2(q-1)} (\wf'+t (\wf-\wf' ) ) \big]\rd t \Big|\nonumber\\
&\lesssim_{q} |\wf-\wf' |\int_0^1 |\wf'+t (\wf-\wf' ) |^{2(q-1)}\rd t\nonumber\\
&\leq ( |\wf |+ |\wf' |)^{2(q-1)}|\wf-\wf' |
\lesssim_{q} (\dens^{q-1}
+\dens'^{q-1} ) |\wf-\wf' |.\nonumber
\end{align*}
Since $H^2$ is embedded into $L^\infty$,
\begin{align}
\label{eq:rhoinfty}
 \|\dens \|_{L^\infty}^{a}\lesssim_{a} \sum_{\elindex=1}^\nel \|\wf_\elindex \|_{H^2}^{2a}
\end{align}
for all $a>0$. Taking $a=q-1>0$ \eqref{eq:D} then follows.\\\\
\textit{Proof of \eqref{eq:E}.}\\
Taking $p=2$ in \eqref{eq:D}, we only need the $\ltwoN$ norm of $\lap(\potX[\dens]\wf-\potX[\dens']\wf' )$ in addition to get the $ \htwoN$ norm estimate. Using the product rule for the Laplacian in $\bb{R}^3$, we get
\begin{align}
\nonumber
\lap(\potX[\dens]\wf-\potX[\dens']\wf' )&=\lambda\Big\{\underbrace{\dens^{q-1}\lap\wf-\dens'^{q-1}\lap\wf'}_{=:\,\,\text{(I)}}+\underbrace{2 [\nab (\dens^{q-1} )\cdot\nab\wf-\nab (\dens'^{q-1} )\cdot\nab\wf' ]}_{=:\,\,\text{(II)}}\\
&\qquad
+\underbrace{\lap(\dens^{q-1} )\wf-\lap(\dens'^{q-1} )\wf'}_{=:\,\,\text{(III)}}\Big\},\label{eq:laplacianP2}
\end{align}
which is in $\bb{C}^\nel$. We discuss the terms one by one.\\\\
\noindent\textit{Term (I).}\\
By adding and subtracting the term $\dens^{q-1}\lap\wf'$ and using \eqref{eq:MVE} with $a=q-1>1$,
\begin{align}
\nonumber
 |\mrm{(I)} |&\leq  |\dens^{q-1} | |\lap\wf-\lap\wf' |+ |\dens^{q-1}-\dens'^{q-1} | |\lap\wf' |\\
&\lesssim_{q} A_1 |\lap\wf' | |\wf-\wf' |
+A_2 |\lap\wf-\lap\wf' |,\label{eq:i}
\end{align}
where
\begin{align*}
A_1= \|\dens \|_{L^{\infty}}^{q-3/2}+ \|\dens' \|_{L^{\infty}}^{q-3/2},\qquad  A_2= \|\dens \|_{L^{\infty}}^{q-1}.
\end{align*}
\textit{Term (II).}\\
By adding and subtracting the term $\nab (\dens^{q-1} )\cdot\nab\wf'$,
\begin{align*}
& |\nab (\dens^{q-1} )\cdot\nab\wf-\nab (\dens'^{q-1} )\cdot\nab\wf' |
\leq\\
&\quad |\nab (\dens^{q-1} ) | |\nab\wf-\nab\wf' |+ |\nab\wf' | |\nab (\dens^{q-1} )-\nab (\dens'^{q-1} ) |.\nonumber
\end{align*}
We get
\begin{align*}
 |\nab (\dens^{q-1} ) |&=(q-1) |\dens^{q-2} | |\nab\dens |\stackrel{\eqref{eq:nablarho}}{\lesssim{}}(q-1) \|\dens \|_{L^{\infty}}^{q-3/2} |\nab\wf |.
\end{align*}
Using this and \eqref{eq:MVEnablarho} with $b=q-1>2$,
\begin{align}
 |\mrm{(II)} |&\lesssim_{q} (B_1 |\nab\wf | |\nab\wf' |+
B_2 |\nab\wf' |^2
 ) |\wf-\wf' |\nonumber\\
&\qquad+ (B_3 |\nab\wf |+B_4 |\nab\wf' | ) |\nab\wf-\nab\wf' |,
\label{eq:ii}
\end{align}
where
\begin{align*}
B_1&= \|\dens \|_{L^{\infty}}^{q-2},
\qquad B_2= \|\dens' \|_{L^{\infty}}^{1/2} ( \|\dens \|_{L^{\infty}}^{q-5/2}+ \|\dens' \|_{L^{\infty}}^{q-5/2} ),\\
B_3&= \|\dens \|_{L^{\infty}}^{q-3/2},\qquad B_4= \|\dens \|_{L^{\infty}}^{q-2} \|\dens' \|_{L^{\infty}}^{1/2}.\nonumber
\end{align*}
\textit{Term (III).}\\
By adding and subtracting the term $\lap(\dens^{q-1} )\wf'$, 
\begin{align*}
 |\mrm{(III)} |\leq\underbrace{ |\lap(\dens^{q-1} ) |}_{=:\,\mrm{(a)}} |\wf-\wf' |+\underbrace{ |\lap(\dens^{q-1} )-\lap(\dens'^{q-1} ) |}_{=:\,\mrm{(b)}} \|\dens' \|^{1/2}_{L^{\infty}}.
\end{align*}
Using $\lap\dens=\overline{\wf}\cdot\lap\wf+\overline{\lap\wf}\cdot\wf+2 |\nab\wf |^2$ and the Cauchy--Schwarz inequality,
\begin{align}
\label{eq:Deltarho}
 |\lap\dens |\lesssim{}  \|\dens \|^{1/2}_{L^{\infty}} |\lap\wf |+ |\nab\wf |^2.
\end{align}
Hence
\begin{align*}
\mrm{(a)}&\lesssim_{q} (q-2) |\dens |^{q-3} |\nab\dens |^2+ |\dens |^{q-2} |\lap\dens |\stackrel{\eqref{eq:nablarho}}{\lesssim{}}
 \|\dens \|_{L^{\infty}}^{q-2}
 \big[
 (2q-3 )
 |\nab\wf |^2
+
 \|\dens \|_{L^{\infty}}^{1/2}
 |\lap\wf |
 \big].
\end{align*}
By similar reasoning, we get, by adding and subtracting the terms $\dens^{q-3} |\nab\dens' |^2$ and~$\dens^{q-2}\lap\dens'$,
\begin{align*}
\mrm{(b)}&\lesssim_{q}(q-2)\big( |\dens |^{q-3}\underbrace{ | |\nab\dens |^2- |\nab\dens' |^2 |}_{=:\,\mrm{(i)}}+\underbrace{ |\dens^{q-3}-\dens'^{q-3} | |\nab\dens' |^2\big)}_{=:\,\mrm{(ii)}}
+
 |\dens |^{q-2}\underbrace{ |\lap\dens-\lap\dens' |}_{=:\,\mrm{(iii)}}\nonumber\\
&\qquad+\underbrace{ |\dens^{q-2}-\dens'^{q-2} | |\lap\dens' |}_{=:\,\mrm{(iv)}}.
\end{align*}
By \eqref{eq:nablarho}  we get
\begin{align*}
\mrm{(i)}&\leq
 (
 |\nab\dens |+ |\nab\dens' |
 )
 |\nab\dens-\nab\dens' |\\
&\stackrel{\eqref{eq:nablarhodiff}}{\lesssim{}} ( \|\dens \|^{1/2}_{L^{\infty}} |\nab\wf |+ \|\dens' \|^{1/2}_{L^{\infty}} |\nab\wf' | ) \big( |\wf-\wf' | |\nab\wf |+
 \|\dens' \|^{1/2}_{L^{\infty}} |\nab\wf-\nab\wf' | \big).\nonumber
\end{align*}
Furthermore, using \eqref{eq:MVE} with $a=q-3\geq1/2$\footnote{Here is where we use the restriction $q\geq 7/2.$ } and \eqref{eq:nablarho} for $\dens'$,
\begin{align*}
\mrm{(ii)}\lesssim_{q}
 ( \|\dens \|^{q-7/2}_{L^{\infty}}+ \|\dens' \|^{q-7/2}_{L^{\infty}} ) \|\dens' \|_{L^{\infty}} |\nab\wf' |^2 |\wf-\wf' |.
\end{align*}
In addition, by adding and subtracting the terms $\overline{\wf'}\cdot\lap\wf$ and $\overline{\lap\wf}\cdot\wf'$, using the triangle and the Cauchy--Schwarz inequalities,
\begin{align*}
\mrm{(iii)}&=\big|2 ( |\nab\wf |^2- |\nab\wf' |^2 )+ (\overline{\wf}-\overline{\wf'} )\cdot\lap\wf+ (\lap\wf-\lap\wf' )\cdot\overline{\wf'}\\
&\quad+ (\wf-\wf' )\cdot\overline{\lap\wf}+ (\overline{\lap\wf}-\overline{\lap\wf'} )\cdot\wf'\big|\nonumber\\
&\lesssim{} ( |\nab\wf |+ |\nab\wf' | ) |\nab\wf-\nab\wf' |+ |\lap\wf | |\wf-\wf' |+ \|\dens' \|_{L^{\infty}}^{1/2} |\lap\wf-\lap\wf' |.\nonumber
\end{align*}
Furthermore, using \eqref{eq:MVE} with $a=q-2>1$ and \eqref{eq:Deltarho} for $\dens'$, we obtain
\begin{align*}
\mrm{(iv)}\lesssim_{q} ( \|\dens \|^{q-5/2}_{L^{\infty}}+ \|\dens' \|^{q-5/2}_{L^{\infty}} ) ( \|\dens' \|_{L^{\infty}}^{1/2} |\lap\wf' |+ |\nab\wf' |^2 ) |\wf-\wf' |.
\end{align*}
Altogether, we get
\begin{align}
\nonumber
 |\mrm{(III)} |&\lesssim_{q}
\big(C_1 |\nab\wf |^2+C_2 |\nab\wf | |\nab\wf' |+C_3 |\nab\wf' |^2+C_4 |\lap\wf |+C_5 |\lap\wf' |
\big) |\wf-\wf' |\nonumber\\
&\quad+ (C_6 |\nab\wf |+C_7 |\nab\wf' | ) |\nab\wf-\nab\wf' |+C_8 |\lap\wf-\lap\wf' |,
\label{eq:iii}
\end{align}
where
\begin{align*}
C_1&= \|\dens \|_{L^{\infty}}^{q-5/2} ( \|\dens \|_{L^{\infty}}^{1/2}+ \|\dens' \|_{L^{\infty}}^{1/2} ),\quad 
C_2= \|\dens \|_{L^{\infty}}^{q-3} \|\dens' \|_{L^{\infty}},\\
C_3&= \|\dens' \|_{L^{\infty}} \big[ \|\dens \|_{L^{\infty}}^{q-7/2} (1+ \|\dens \|_{L^{\infty}} )+ \|\dens' \|_{L^{\infty}}^{q-7/2} (1+ \|\dens' \|_{L^{\infty}} ) \big],\\
C_4&= \|\dens \|_{L^{\infty}}^{q-3/2} \big( \|\dens \|_{L^{\infty}}^{1/2} \|\dens' \|_{L^{\infty}}^{1/2}+1 \big),\quad C_5= \|\dens' \|_{L^{\infty}} \big( \|\dens \|_{L^{\infty}}^{q-5/2}+ \|\dens' \|_{L^{\infty}}^{q-5/2} \big),\\
C_6&= \|\dens \|_{L^{\infty}}^{q-5/2} \|\dens' \|_{L^{\infty}}^{1/2} ( \|\dens \|_{L^{\infty}}^{1/2}+ \|\dens' \|_{L^{\infty}}^{1/2} ),\quad 
C_7= \|\dens \|_{L^{\infty}}^{q-3} \|\dens' \|_{L^{\infty}}^{1/2} ( \|\dens \|_{L^{\infty}}+ \|\dens' \|_{L^{\infty}} ),\\
C_8&= \|\dens \|_{L^{\infty}}^{q-2} \|\dens' \|_{L^{\infty}}.
\end{align*}
\textit{Conclusion of the proof of \eqref{eq:E}.}\\
The function $\mathcal{L}$ can be split into terms $\mathcal{L}=\mathcal{L}_0+\mathcal{L}_{\mrm{I}}+\mathcal{L}_{\mrm{II}}+\mathcal{L}_{\mrm{III}}$. As discussed at the start of this proof, $\mathcal{L}_0$ is the contribution of estimate \eqref{eq:D} for $p=2$. The other terms stem from (I), (II) and (III) in \eqref{eq:laplacianP2}, and are obtained taking the $L^2$ norm in (\ref{eq:i}, \ref{eq:ii} resp. \ref{eq:iii}). For instance, in the expression of $\mathcal{L}_{\mrm{III}}$ all scalars $C_i$ can be bounded using \eqref{eq:rhoinfty}. Likewise, the term  $|\wf-\wf' |$  as well as the remaining factors involving $C_1$ and $C_3$ can be bounded with their $H^2$ norms. The other summands, and $\mathcal{L}_{\mrm{I}}$ and $\mathcal{L}_{\mrm{II}},$ can be estimated similarly, and this concludes the proof.
\end{proof}
\begin{lemma}[Lipschitz estimates for the nonlinearity]
\label{lem:lipschitzP}
For $q\geq 7/2$ and any $\lambda\in\bb{R}$, there exists a polynomial function which vanishes at the origin $\mathscr{L}_{q,\lambda}:\bb{R}_0^+\longrightarrow \bb{R}_0^+$, such that for all $\wf,\wf'\in\Bel$
\begin{align}
&\|\pot_{\text{Hx}}[\dens]\wf-\pot_{\text{Hx}}[\dens']\wf' \|_{C^0 ([0,\tau]; \htwoN )}\leq\nonumber\\
&\quad\leq \mathscr{L}_{q,\lambda} (\alpha(\tau)+ \|\wf^0 \|_{ \htwoN} ) \|\wf-\wf' \|_{C^0 ([0,\tau]; \htwoN )},\label{eq:F}\\
&\|\pot_{\text{Hx}}[\dens]\wf \|_{C^0 ([0,\tau]; \htwoN )}\leq\nonumber\\
&\quad \leq (\alpha(\tau)+ \|\wf^0 \|_{ \htwoN} )\mathscr{L}_{q,\lambda} (\alpha(\tau)+ \|\wf^0 \|_{ \htwoN} ).\label{eq:G}
\end{align}
\end{lemma}
\begin{proof}
By \eqref{eq:C} in Lemma \ref{lem:lipschitzconvolution} and \eqref{eq:E} in Lemma \ref{lem:lipschitzlda}, we have, for all $\wf,\wf'\in C^0 ([0,\tau]; \htwoN ),$ that 
\begin{align}
&\|\pot_{\text{Hx}}[\dens]\wf-\pot_{\text{Hx}}[\dens']\wf' \|_{C^0 ([0,\tau]; \htwoN )}\leq\nonumber\\
&\leq \mathscr{L}_{q,\lambda} (\max \{ \|\wf \|_{C^0 ([0,\tau]; \htwoN )}, \|\wf' \|_{C^0 ([0,\tau]; \htwoN )} \} )\|\wf-\wf' \|_{C^0 ([0,\tau]; \htwoN )},\label{eq:LipP}
\end{align}
where $\mathscr{L}_{q,\lambda}$ is a polynomial by construction. Note that \eqref{eq:F} follows from \eqref{eq:LipP} by definition of $\Bel$. In particular,  \eqref{eq:G} follows from \eqref{eq:LipP} setting $\wf'\equiv 0$.
\end{proof}
\section{Existence and uniqueness of nuclear configurations}
\label{sec:locex}
In this section, we prove a local-in-time existence and uniqueness result for the Cauchy problem associated with \eqref{eq:N} for given $\wf\in\Bel$ and $\nuc(0)=\nuc^0,\dot{\nuc}(0)=\vel(0)$, with $\nuc^0,\vel^0\in\bb{R}^{3\nnuc}$ such that $\nuc^0_\nucindex\neq\nuc^0_\nucindextwo$ for $1\leq\nucindex\neq\nucindextwo\leq\nnuc$.
\begin{lemma}
\label{lem:KStoCL}
Let $\vel^0\in\bb{R}^{3\nnuc}$ and $\nuc^0\in\bb{R}^{3\nnuc}$ be given, with $\nuc^0_{\nucindex}\neq \nuc^0_{\nucindextwo}$ for $\nucindex\neq \nucindextwo$.\newline 
Then, there exists $\tau>0$ such that the following properties hold. For given $\wf\in\Bel$, the Cauchy problem associated with the system \eqref{eq:N} with $\nuc(0)=\nuc^0$ and $\dot{\nuc}(0)=\vel^0$ has a unique short-time solution $\nuc\in\Bnuc\cap C^2 ([0,\tau];B_\delta (\nuc^0 ) )$. The mapping
\begin{align*}
\mathcal{N}:\wf \in \mathcal{B}_{\mathrm{el}} (\tau )\longmapsto \nuc\in\mathcal{B}_{\mathrm{nuc}} (\tau )\cap C^2 ([0,\tau];\bb{R}^{3\nnuc} ) 
\end{align*}
is bounded with respect to the $C^1 ([0,\tau];\bb{R}^{3\nnuc} )$ topology, and continuous as a map from $C^0\big([0,\tau];\ltwoN\big)$ to $C^0\big([0,\tau];\bb{R}^{3\nnuc}\big)$.
\end{lemma}
\begin{proof}
\textit{Part 1: Existence and uniqueness of $\nuc$ in $C^2 ([0,\tau];B_\delta (\nuc^0 ) )$.}\\
Since~$\wf$ and so $\dens$ are given, we write the acceleration function from \eqref{eq:N} without parameters for now: $\acc=\acc(t,\nuc)$. Note that $t$ is an explicit variable for the $\acc^1_{\nucindex}$ terms, but not for the $\acc^2_{\nucindex}$ terms.\newline
We define the compact set 
\begin{align*}
\varkappa(\tau):=[0,\tau]\times B_\delta(\nuc^0). 
\end{align*}
Note we drop the dependence of this set on $\tau$. By the triangle inequality, for all~$\nuc\in B_{\delta} (\nuc^0 )$ and $\nucindex=1,\ldots,\nnuc$
\begin{align}
 | \nuc_{\nucindex} |&
\leq |\nuc |\leq  |\nuc^0 |+ |\nuc-\nuc^0 |\leq  |\nuc^0 |+\delta(\tau).
\label{eq:delta2}
\end{align}
First, we prove $\acc$ is continuous in $(t,\nuc)$ on $\varkappa$. To this end, we pick a sequence $ \{(t_n,\nuc_n) \}_{n\in\bb{N}}\subset \varkappa$ with~$(t_n,\nuc_n)\xrightarrow{n\longrightarrow\infty}(t^*,\nuc)\in \varkappa$. The functions $\acc^1_{\nucindex}$ give for all $n$, using the Cauchy--Schwarz and Hardy's inequalities,
\begin{align*}
&|\acc^1_{\nucindex} (t_n,\nuc_n )-\acc^1_{\nucindex} (t^*,\nuc_n ) |\lesssim{}\nonumber\\
&\quad \lesssim{}
\sum_{\elindex=1}^\nel \langle |\nuc_{n\nucindex}-\cdot\, |^{-2}, | (\wf_\elindex(t_n,\cdot) )^2- (\wf_\elindex(t^*,\cdot) )^2 | \rangle_{L^2}\nonumber\\
&\quad \lesssim{} \sum_{\elindex=1}^\nel\big\langle |\nuc_{n\nucindex}-\cdot\, |^{-1}\max_{t\in[0,\tau]} |\wf_\elindex(t,\cdot) | ,  |\nuc_{n\nucindex}-\cdot\, |^{-1} |\wf_\elindex(t_n,\cdot)-\wf_\elindex(t^*,\cdot) |\big\rangle_{L^2}\nonumber\\
&\quad \lesssim{}\sum_{\elindex=1}^\nel \|\nab\wf_\elindex \|_{L^{\infty} ([0,\tau];L^2 (\bb{R}^3;\bb{C}^3 ) )} \|\nab\wf_\elindex(t_n,\cdot)-\nab\wf_\elindex(t^*,\cdot) \|_{L^2 (\bb{R}^3;\bb{C}^3 )}\xrightarrow{n\longrightarrow\infty}0,
\end{align*}
as $\wf\in C^0 ([0,\tau]; \honeN )$. Using this and Lemma \ref{lem:forces}, by which $\acc^1_{\nucindex}(t^*,\cdot)\in C^0 (\bb{R}^{3\nnuc};\bb{C}^3 )$, 
\begin{align*}
 |\acc^1_{\nucindex}(t_n,\nuc_n)-\acc^1_{\nucindex}(t^*,\nuc) |
&\leq  |\acc^1_{\nucindex}(t_n,\nuc_n) -\acc^1_{\nucindex}(t^*,\nuc_n) | +  |\acc^1_{\nucindex}(t^*,\nuc_n)-\acc^1_{\nucindex}(t^*,\nuc) |\nonumber\\
&\qquad\xrightarrow{n\longrightarrow\infty}0
\end{align*}
for all $n$. The functions $\acc^2_{\nucindex}$ are not explicitly time-dependent, and are continuous on~$B_\delta (\nuc^0 )$, hence on $\varkappa$. 

Since $\acc$ is continuous on the compact set $\varkappa$, it is also uniformly bounded on $\varkappa$.
By Lemma \ref{lem:forces},
\begin{align*}
 \|\acc^1_{\nucindex} \|_{C^0 ([0,\tau];W^{1,\infty} (B_\delta (\nuc^0 );\bb{C}^3 ) )}&\lesssim{} \|\wf \|^2_{C^0 ([0,\tau]; \htwoN )},
\end{align*}
since $\wf\in\Bel$. The functions $\acc^2_\nucindex$ are bounded on $B_\delta (\nuc^0 )$ by
\begin{align*}
 \|\acc^2_\nucindex \|_{L^{\infty} (B_\delta (\nuc^0 );\bb{C}^3 )}
&\lesssim{}
\sum_{\nucindextwo=1,\nucindextwo\neq \nucindex}^{\nnuc} \Big\|\frac{1}{ | \nuc_{\nucindex}-\nuc_{\nucindextwo} |^{2}} \Big\|_{L^{\infty} (B_\delta (\nuc^0 );\bb{C} )}.
\end{align*}
Furthermore, by Lemma \ref{lem:forces}, $\acc^1_{\nucindex}(t,\cdot)$ is uniformly Lipschitz continuous for all $t\in[0,\tau]$ and $\nucindex$, as
\begin{align*}
 \|D\acc_{\nucindex}^1(t,\cdot) \|_{L^{\infty}(\bb{R}^3;\bb{C}^{3\times 3})}&\lesssim \sum_{\elindex=1}^\nel \|Df^{\elindex\elindex}_\nucindex(t,\cdot) \|_{L^{\infty}(\bb{R}^3;\bb{C}^{3\times 3})}\lesssim{} \|\wf \|^2_{C^0 ([0,\tau]; \htwoN )}
\end{align*}
since $\wf\in\Bel$. For the $\acc^2_{\nucindex}$ terms, we note that the functions $\nuc\longmapsto (\nuc_\nucindex-\nuc_{\nucindextwo} )|\nuc_\nucindex-\nuc_{\nucindextwo}|^{-3}$ are locally Lipschitz on $B_{\delta} (\nuc^0 )$. Therefore, $\acc$ is Lipschitz continuous in $\nuc\in B_{\delta} (\nuc^0 )$ and uniformly in $t\in[0,\tau]$. We denote the corresponding Lipschitz constant by $C_\text{L}$, dropping its dependence on $\tau.$ 

Now, we define $\mathcal{T}$ as the following mapping on the complete metric space $C^0 ([0,\tau];B_\delta (\nuc^0 ) )$, equipped with the $C^0 ([0,\tau];\bb{R}^{3\nnuc} )$ norm:
\begin{align}
\label{eq:mappingT}
\mathcal{T} [\nuc ](t):=\nuc^0+\vel^0t+\int_0^t(t-\sigma)\acc(\sigma,\nuc(\sigma))\rd \sigma.
\end{align}
By the boundedness of $\acc$,
\begin{align*}
 \|\mathcal{T} [\nuc ]-\nuc^0 \|_{C^0 ([0,\tau];\bb{R}^{3\nnuc} )}\leq  |\vel^0 |\tau+\frac{\tau^2}{2} \|\acc \|_{C^0 (\varkappa;\bb{C}^{3\nnuc} )}
\end{align*}
for all $\nuc\in C^0 ([0,\tau];B_\delta (\nuc^0 ) )$. Note that $\mathcal{T}$ maps $C^0 ([0,\tau];B_\delta (\nuc^0 ) )$
into itself, as for $\tau>0$ small enough it holds that
\begin{align*}
 |\vel^0 |\tau+
\frac{\tau^2}{2} \|\acc \|_{C^0 (\varkappa;\bb{C}^{3\nnuc} )}\leq\delta(\tau).
\end{align*}
Hence, for all $\nuc,\nuc'\in C^0([0,\tau];B_{\delta} (\nuc^0 ))$,
\begin{align*}
 \|\mathcal{T} [\nuc ]-\mathcal{T} [\nuc' ] \|_{C^0 ([0,\tau];\bb{R}^{3\nnuc} )}&\leq \max_{t\in[0,\tau]} \int_0^t(t-\sigma) |\acc (\sigma,\nuc(\sigma) )-\acc (\sigma,\nuc'(\sigma) ) |\rd \sigma \nonumber\\
&\leq\frac{C_\text{L}\tau^2}{2} \|\nuc-\nuc' \|_{C^0 ([0,\tau];\bb{R}^{3\nnuc} )},
\end{align*}
Note also that $\mathcal{T}$ is a strict contraction on $C^0 ([0,\tau];B_\delta (\nuc^0 ) )$ in the $C^0 ([0,\tau];\bb{R}^{3\nnuc} )$ norm, as we can always shrink $\tau>0$ so that
\begin{align*}
\frac{C_\text{L}\tau^2}{2} < 1
\end{align*}
holds. By the contraction mapping theorem, $\mathcal{T}$ has a unique fixed point in $C^0([0,\tau];B_{\delta} (\nuc^0 ))$. Because of this, \eqref{eq:N} has a unique short-time solution in $C^2 ([0,\tau];B_{\delta} (\nuc^0 ) )$.\\\\
\textit{Part 2: Localisation of $\nuc$ in $\Bnuc$.}\\
Integrating the ODE in \eqref{eq:N}, we get
\begin{align*}
 \|\dot{\nuc} \|_{C^0 ([0,\tau];\mathbb{R}^{3\nnuc} )}
&\leq 
 |\vel^0 |+\tau \|\acc \|_{C^0 (\varkappa;\bb{C}^{3\nnuc} )}.
\end{align*}
Note that $\nuc\in \Bnuc$, picking $\tau>0$ smaller if necessary, so that
\begin{align*}
\tau \|\acc \|_{C^0 (\varkappa;\bb{C}^{3\nnuc} )}\leq 1,
\end{align*}
holds. Therefore, $\nuc\in \Bnuc\cap C^2 ([0,\tau];B_\delta (\nuc^0 ) )$.\\\\
\textit{Part 3: Boundedness and continuity of $\mathcal{N}$.}\\
From \eqref{eq:delta2} with \eqref{eq:gamma}, $\mathcal{N}$ is bounded in the $C^1 ([0,\tau];\bb{R}^{3\nnuc} )$ norm:
\begin{align*}
 \|\nuc \|_{C^1 ([0,\tau];\bb{R}^{3\nnuc} )}
&\leq  |\nuc^0 |+\delta(\tau)+\gamma.
\end{align*}
In order to prove continuity of $\mathcal{N}$ in the $C^0 ([0,\tau];\bb{R}^{3\nnuc} )$ norm, we consider a sequence~$ \{\wf_n \}_{n\in\bb{N}}\subset\Bel$ such that $\wf_n\xrightarrow{n\longrightarrow\infty}\wf\in\Bel$ in the $C^0 ([0,\tau];\ltwoN )$ norm. Similarly to $\nuc=\mathcal{N} [\wf ]$, we define $\nuc_n:=\mathcal{N} [\wf_n ]$ and $\dens_n:= |\wf_n |^2$. Note that $\nuc$ and $\nuc_n$ are fixed points of the mapping $\mathcal{T}$ introduced in Part 1 of the proof. Using this, for all $t\in[0,\tau]$
\begin{align}
& | (\nuc_{n}-\nuc )(t) |
\leq
\int_0^t(t-\sigma) |\acc [\dens_n ](\nuc(\sigma))-\acc [\dens ](\nuc(\sigma)) |\rd\sigma,
\label{eq:C0cty}
\end{align}
where
\begin{align*}
 |\acc [\dens_n ] (\nuc_n(\sigma) )-\acc [\dens ](\nuc(\sigma)) |&\leq
\mrm{(I)}+\mrm{(II)},\\
\mrm{(I)}&:=\sum_{\nucindex=1}^{\nnuc} |\acc^1_{\nucindex} [\dens_n ] (\nuc_{n}(\sigma) )-\acc^1_{\nucindex} [\dens ](\nuc(\sigma)) |,\nonumber\\
\mrm{(II)}&:=\sum_{\nucindex=1}^{\nnuc} |\acc^2_{\nucindex} (\nuc_{n}(\sigma) )-\acc^2_{\nucindex} (\nuc(\sigma) ) |.\nonumber
\end{align*}
We further bound
\begin{align*}
\mrm{(I)}&\lesssim{}\mrm{(Ia)}+\mrm{(Ib)},\\
\mrm{(Ia)}&:=\sum_{\elindex=1}^\nel\Big| \langle\wf_{\elindex}(t,\cdot),\Xi(\cdot- \nuc_{\nucindex}) (\wf_{n\elindex}(t,\cdot)-\wf_\elindex(t,\cdot) ) \rangle_{L^2}\nonumber\\&\qquad + \langle\wf_{n\elindex}(t,\cdot)-\wf_\elindex(t,\cdot),\Xi(\cdot-\nuc_{n\nucindex} )\wf_{n\elindex}(t,\cdot) \rangle_{L^2}\Big|,\nonumber\\
\mrm{(Ib)}&:=\sum_{\elindex=1}^\nel | \langle\wf_\elindex(t,\cdot),\Xi(\cdot-\nuc_{n\nucindex} )\wf_{n\elindex}(t,\cdot) \rangle_{L^2}- \langle\wf_{\elindex}(t,\cdot),\Xi(\cdot- \nuc_{\nucindex} )\wf_{n\elindex}(t,\cdot) \rangle_{L^2} |,\nonumber
\end{align*}
where we use for short-hand notation the function $\Xi:\bb{R}^3\longrightarrow \bb{R}^3$ (a.e.),~$\el\longmapsto\el |\el |^{-3}$. Arguing as in \cite[p. 980]{Cances1999OnDynamics}, (Ia) is bounded by
\begin{align*}
\beta_{n}&:=\sum_{\elindex=1}^\nel
\sup_{(t,\el)\in[0,\tau]\times\bb{R}^3}
 \langle
 |\cdot-\el |^{-1} |\wf_\elindex(t,\cdot)+\wf_{n\elindex}(t,\cdot) |,
 |\cdot-\el |^{-1} |\wf_{n\elindex}(t,\cdot)-\wf_\elindex(t,\cdot) | \rangle_{L^2}
 \nonumber\\
&\qquad \xrightarrow{n\longrightarrow\infty}0,
\end{align*}
as $\wf_n\xrightarrow{n\longrightarrow\infty}\wf$ in $C^0 ([0,\tau];\ltwoN )$. We also have
\begin{align*}
\mrm{(Ib)}&\lesssim{}\sum_{\elindex=1}^\nel \|\nab G [\wf_\elindex,\wf_{n\elindex} ] (\nuc_{n\nucindex} )-
\nab G [\wf_\elindex,\wf_{n\elindex} ] ( \nuc_{\nucindex} ) \|_{C^0 ([0,\tau];\bb{C}^3 )}\leq C^\text{L}_{1,n} |\nuc_{n}-\nuc |,
\end{align*}
where $G$ is as in \eqref{eq:mappingG}, and where we used that the functions $\nab G [\wf_\elindex,\wf_{n\elindex} ]$ are uniformly Lipschitz continuous in $\nuc$ for uniformly all $t\in[0,\tau]$. So is (II), with some Lipschitz constant $C^\text{L}_{2,n}$. For all $n$, $C^\text{L}_{1,n}$ and $C^\text{L}_{2,n}$ are uniformly bounded by $C_\text{L}$, since all $\wf_n$ and $\wf$ are taken from the bounded set $\Bel$. Altogether, from \eqref{eq:C0cty} we obtain
\begin{align*}
 \|\nuc_{n}-\nuc \|_{C^0 ([0,\tau];\bb{R}^{3\nnuc} )}\lesssim \tau^2 \|\nuc_{n}-\nuc \|_{C^0 ([0,\tau];\bb{R}^{3\nnuc} )}+\tau^2\beta_n.
\end{align*}
It is then clear that for $\tau$ small enough the conclusion follows.
\qedhere
\end{proof}
\section{Existence and uniqueness of electronic configurations}
\label{sec:el}
In this section, we prove a local-in-time existence and uniqueness result for the Cauchy problem associated with \eqref{eq:KS} for given $\nuc\in\Bnuc$ and $\wf(0)=\wf^0\in \htwoN$.
\begin{lemma}
\label{lem:CLtoKS}
Let $q\geq 7/2$, $\lambda\in\bb{R}$. Let $\wf^0\in \htwoN$ be given. Then, there exists $\tau>0$ such that the following holds. For given $\nuc\in\Bnuc$, the Cauchy problem associated with the system \eqref{eq:KS} with $\wf(0)=\wf^0$ has a unique short-time solution $\wf$ in $\Bel$.
\end{lemma}
\begin{proof}
This proof is based on Lemma \ref{lem:evolutionoperators}, which ensures the existence and the $\mathcal{L} (\htwoN )$ bounds of the propagator $U(t,s)$ for the family of linear Hamiltonians $ \{\schoplin(t),t\in[0,\tau] \}$ from \eqref{eq:linHam}, and on Lemma \ref{lem:lipschitzP}, which ensures that the nonlinear mapping $\wf\longmapsto \pot_{\text{Hx}}[\dens]\wf$ is locally Lipschitz in $ \htwoN$.\newline
We define $\mathcal{F}$ as the following mapping on the complete metric space $C^0 ([0,\tau];B_{\alpha} (\wf^0 ) )$, equipped with the $C^0 ([0,\tau]; \htwoN )$ norm:
\begin{align*}
\mathcal{F}[\wf] (t):=U(t,0)\wf^0-i\int_0^tU(t,\sigma)\pot_{\text{Hx}}[\dens]\wf(\sigma)\rd\sigma.
\end{align*}
Note that we obtain for all $\wf\in C^0 ([0,\tau];B_{\alpha} (\wf^0 ) )$, using Lemma \ref{lem:evolutionoperators} (ii),
\begin{align}
\mathcal{F}[\wf](0)&=U(0,0)\wf^0=\wf^0.\label{eq:fpsi0}
\end{align}
Note also that, provided 
\begin{align}
\big[1+B_{\tau,\gamma}+\tau B_{\tau,\gamma} (2B_{\tau,\gamma}+1 )\mathscr{L}_{q,\lambda} (\alpha+ \|\wf^0 \|_{ \htwoN} ) \big]\leq2B_{\tau,\gamma},\label{eq:A2a}
\end{align}
$\mathcal{F}$ maps the complete metric space $C^0 ([0,\tau];B_{\alpha} (\wf^0 ) )$ into itself, as
\begin{align}
&\|\mathcal{F}[\wf]-\wf^0 \|_{C^0 ([0,\tau]; \htwoN )}
=\nonumber\\
&\quad =\Big\| [U(\cdot,0)-\mrm{Id} ]\wf^0-i\int_0^\cdot U(\cdot,\sigma)\pot_{\text{Hx}}[\dens]\wf(\sigma)\rd\sigma \Big\|_{C^0 ([0,\tau]; \htwoN )}\nonumber\\
&\quad \leq B_{\tau,\gamma} \big( \|\wf^0 \|_{ \htwoN}+\tau  \|\pot_{\text{Hx}}[\dens]\wf \|_{C^0 ([0,\tau]; \htwoN )} \big) + \|\wf^0 \|_{ \htwoN}\label{eq:star14}\\
&\quad \stackrel{\eqref{eq:G}}{\leq}\big[1+B_{\tau,\gamma}+\tau B_{\tau,\gamma} (2B_{\tau,\gamma}+1 )\mathscr{L}_{q,\lambda} (\alpha+ \|\wf^0 \|_{ \htwoN} ) \big]  \|\wf^0 \|_{ \htwoN}\nonumber\\
&\quad \stackrel{\eqref{eq:A2a}}{\leq} 2B_{\tau,\gamma} \|\wf^0 \|_{ \htwoN}=\alpha,\nonumber
\end{align}
where we used Lemma \ref{lem:evolutionoperators} (viii) in \eqref{eq:star14}. 
Moreover, note that, provided 
\begin{align}
\tau B_{\tau,\gamma}\mathscr{L}_{q,\lambda} (\alpha(\tau)+ \|\wf^0 \|_{ \htwoN} )
<1,\label{eq:A2b}
\end{align}
$\mathcal{F}$ is a contraction on $C^0 ([0,\tau];B_{\alpha} (\wf^0 ) )$ in the~$C^0 ([0,\tau]; \htwoN )$ norm, as for all $\wf,\wf'\in C^0 ([0,\tau];B_{\alpha} (\wf^0 ) )$
\begin{align}
&\|\mathcal{F}[\wf]-\mathcal{F}[\wf'] \|_{C^0 ([0,\tau]; \htwoN )}=\nonumber\\
&\quad = \bigg\|\int_0^\cdot U(\cdot,\sigma) (\pot_{\text{Hx}}[\dens]\wf(\sigma)-\pot_{\text{Hx}}[\dens']\wf'(\sigma) )\rd\sigma \bigg\|_{C^0 ([0,\tau]; \htwoN )}\nonumber\\
&\quad \leq\tau B_{\tau,\gamma} \|\pot_{\text{Hx}}[\dens]\wf-\pot_{\text{Hx}} [\dens' ]\wf' \|_{C^0 ([0,\tau]; \htwoN )}\label{eq:star15}\\
&\quad \stackrel{\eqref{eq:F}}{\leq}\tau B_{\tau,\gamma}\mathscr{L}_{q,\lambda} (\alpha+ \|\wf^0 \|_{ \htwoN} ) \|\wf-\wf' \|_{C^0 ([0,\tau]; \htwoN )},\nonumber
\end{align}
where we used Lemma \ref{lem:evolutionoperators} (viii) in \eqref{eq:star15}. By the contraction mapping theorem, $\mc{F}$ has a unique fixed point in $C^0 ([0,\tau];B_{\alpha} (\wf^0 ) )$.

Note that we can always select $\tau>0$ small enough such that the inequalities (\ref{eq:A2a},\ref{eq:A2b}) are satisfied. Recall that $B_{\tau,\gamma}$ and $\alpha$ are of the form
\begin{align*}
B_{\tau,\gamma}=A_\gamma^{1+C_\gamma\tau},\qquad \alpha(\tau)=2A_\gamma^{1+C_\gamma\tau} \|\wf^0 \|_{ \htwoN},
\end{align*}
with $A_\gamma,C_\gamma>0$ as defined as in Lemma \ref{lem:evolutionoperators} (viii). In fact, picking $A_\gamma,C_\gamma$ large, (\ref{eq:A2a}) is true for $\tau=0$ and by continuity, for $\tau>0$ small enough. \\

It is now left to prove that this fixed point, simply denoted by $\wf$, is also of class~$C^1 ([0,\tau];\ltwoN )$; namely, it solves \eqref{eq:KS} strongly on $[0,\tau]$. To this end, we consider the following identity, which holds for all $0\leq t<t'\leq \tau$:
\begin{align*}
i\frac{\wf(t')-\wf(t)}{t'-t}=\mrm{(R)}&:=
i\frac{U(t',0)-U(t,0)}{t'-t}\wf^0+\int_0^t\frac{U(t',\sigma)-U(t,\sigma)}{t'-t}\pot_{\text{Hx}} [\dens ]\wf(\sigma)\rd\sigma\nonumber
\\
&\qquad +\int_t^{t'}\frac{U(t',\sigma)}{t'-t}\pot_{\text{Hx}} [\dens ]\wf(\sigma)\rd\sigma,
\end{align*}
and we show that 
\begin{align*}
 \|\mrm{(R)}-\schop [\nuc(t),\dens ]\wf(t) \|_{\ltwoN}\xrightarrow{t'\longrightarrow t}0.
\end{align*}
This will imply that $ \wf(\cdot)$ is differentiable as a mapping $[0,\tau] \longmapsto \ltwoN$  such that 
\begin{align*}
i\dot{\wf}(t)=\schop [\nuc(t),\dens ]\wf(t).
\end{align*}
Note, in particular, that for a given $\nuc\in\Bnuc$, $\schop [\nuc(\cdot),\dens ]\wf(\cdot)$ is a continuous mapping~$[0,\tau] \longmapsto \ltwoN$, which will imply that $ \wf\in C^1 ([0,\tau];\ltwoN ).$ We bound 
\begin{align*}
& \|\mrm{(R)}-\schop [\nuc(t),\dens ]\wf(t) \|_{\ltwoN}\leq
\mrm{(I)}+\mrm{(II)},\\
&\mrm{(I)}:= \Big\|
i\frac{U(t',0)-U(t,0)}{t'-t}\wf^0+
\int_0^t\frac{U(t',\sigma)-U(t,\sigma)}{t'-t}\pot_{\text{Hx}} [\dens ]\wf(\sigma)\rd\sigma
-
\schoplin(t)\wf(t)
 \Big\|_{\ltwoN},\nonumber\\
&\mrm{(II)}:= \Big\|
\int_t^{t'}\frac{U(t',\sigma)}{t'-t}\pot_{\text{Hx}} [\dens ]\wf(\sigma)\rd\sigma
-
\pot_{\text{Hx}} [\dens ]\wf(t)
 \Big\|_{\ltwoN}.\nonumber
\end{align*}
We get 
\begin{align}
\lim_{t'\longrightarrow t} \mrm{(I)} &= \Big\|i\frac{\partial}{\partial t} [U(t,0)\wf^0 ]
+\int_0^t
\frac{\partial}{\partial t} [U(t,\sigma)\pot_{\text{Hx}} [\dens ]\wf(\sigma) ]\rd\sigma-
\schoplin(t)\wf(t)
 \Big\|_{\ltwoN}
\nonumber\\
&=
\Big\|
\schoplin(t) [U(t,0)\wf^0 ]\nonumber\\
&\qquad 
+
\int_0^t
-i
\schoplin(t) [
U(t,\sigma)\pot_{\text{Hx}} [\dens ]\wf(\sigma)
 ]
\rd\sigma
-
\schoplin(t)\wf(t)
\biggr\|_{\ltwoN}
\label{eq:star16}\\
&=
 \|
\schoplin(t) [\mathcal{F} [\wf(t) ]-\wf(t) ]
 \|_{\ltwoN}=0,\label{eq:star26}
\end{align}
where we used Lemma \ref{lem:evolutionoperators} (vii) (see also \cite[Thm. 1.3. (6)]{Yajima1987ExistenceEquations}) in \eqref{eq:star16}, the linearity of the Hamiltonians $\schoplin(t)$ and $\wf$ being a fixed point of $\mathcal{F}$ in \eqref{eq:star26}. On the other hand,
\begin{align*}
\mrm{(II)}&\leq
\mrm{(a)}+\mrm{(b)},\\
\mrm{(a)}&:= \Big\|
\frac{1}{t'-t}\int_t^{t'}U(t,\sigma)\pot_{\text{Hx}} [\dens ]\wf(\sigma)\rd\sigma
-
\pot_{\text{Hx}} [\dens ]\wf(t)
 \Big\|_{\ltwoN},\nonumber\\
\mrm{(b)}&:=\frac{1}{t'-t} \Big\|
\int_t^{t'} [U(t',\sigma)-U(t,\sigma) ]\pot_{\text{Hx}} [\dens ]\wf(\sigma)\rd\sigma   
 \Big\|_{\ltwoN}.\nonumber
\end{align*}
In the limit, (a) goes to zero, because of the fundamental theorem of calculus for Bochner integrals and Lemma \ref{lem:evolutionoperators} (ii). Moreover,
\begin{align}
\lim_{t'\longrightarrow t} \mrm{(b)} &\leq\lim_{t'\longrightarrow t} \frac{1}{t'-t}
\int_t^{t'} \| [U(t',\sigma)-U(t,\sigma) ]\pot_{\text{Hx}} [\dens ]\wf(\sigma) \|_{\ltwoN}\rd\sigma \nonumber\\
&\leq\lim_{t'\longrightarrow t}  \| [U(t',\cdot)-U(t,\cdot) ]\pot_{\text{Hx}} [\dens ]\wf \|_{C^0 ([0,T],\ltwoN )} =0,\label{eq:star17}
\end{align}
where we used the uniform continuity of 
$U(t,s)\pot_{\text{Hx}} [\dens ]\wf(s)$ on $[0,T]^2$ together with Lemma \ref{lem:evolutionoperators} (iv) in \eqref{eq:star17}. Since $\wf$ also is a fixed point of $\mathcal{F}$, by which $\wf(0)=\mathcal{F}[\wf](0)=\wf^0$ (see \eqref{eq:fpsi0}), we know $\wf$ is a strong solution to \eqref{eq:KS} on $[0,\tau]$.

Now, we show uniqueness of the short-time solution $\wf$ to \eqref{eq:KS} in the class~$C^1 ([0,\tau];\ltwoN )\cap C^0 ([0,\tau];B_{\alpha} (\wf^0 ) )$: although the classical contraction mapping theorem also provides uniqueness, this is only in the class $C^0 ([0,\tau];B_{\alpha} (\wf^0 ) )$. So, now we prove uniqueness in the different space $C^1 ([0,\tau];\ltwoN )$. To this end, we let~$\wf$ and $\wf'$ be two short-time solutions to \eqref{eq:KS} in $C^1 ([0,\tau];\ltwoN )$. First, $ (\wf-\wf' )(0)=\wf^0-\wf^0=0$. Moreover, for all $\elindex = 1,\ldots,\nel$, using the PDE in \eqref{eq:KS},
\begin{align*}
\frac{\rd}{\rd t} ( \|\wf_{\elindex}-\wf_{\elindex}' \|^2_{L^2} )&=
\frac{\rd}{\rd t} ( \langle \wf_{\elindex}-\wf_{\elindex}',\wf_{\elindex}-\wf_{\elindex}' \rangle_{L^2} )\\
&= \langle \dot{\wf}_{\elindex}-\dot{\wf}_{\elindex}',\wf_{\elindex}-\wf_{\elindex}' \rangle_{L^2}+\overline{ \langle \dot{\wf}_{\elindex}-\dot{\wf}_{\elindex}',\wf_{\elindex}-\wf_{\elindex}' \rangle_{L^2}}=\mrm{(I)}+\mrm{(II)},\nonumber
\end{align*}
where, using that the linear Hamiltonians $\schoplin(t)$ are self-adjoint on $\ltwoN$,
\begin{align*}
\mrm{(I)}&=
i \big[\langle\wf_\elindex-\wf_\elindex', (\schoplin(t) (\wf-\wf' ) )_\elindex \rangle_{L^2}-
 \langle (\schoplin(t) (\wf-\wf' ) )_\elindex,\wf_{\elindex}-\wf_{\elindex}' \rangle_{L^2}
 \big]=0,
\end{align*}
and
\begin{align*}
\mrm{(II)}&=i \big[
\overline{ \langle (\pot_{\text{Hx}} [\dens ]\wf-\pot_{\text{Hx}} [\dens' ]\wf' )_\elindex,\wf_\elindex-\wf_\elindex' \rangle_{L^2}}-
 \langle (\pot_{\text{Hx}} [\dens ]\wf-\pot_{\text{Hx}} [\dens' ]\wf' )_\elindex,\wf_\elindex-\wf_\elindex' \rangle_{L^2} \big]\nonumber
\\
&=2\mathrm{Im} \langle (\pot_{\text{Hx}} [\dens ]\wf-\pot_{\text{Hx}} [\dens' ]\wf' )_\elindex,\wf_\elindex-\wf_\elindex' \rangle_{L^2}.
\end{align*}
Using this, we get
\begin{align*}
\frac{\rd}{\rd t} \big( \|\wf-\wf' \|^2_{\ltwoN} \big)
&=
\sum_{\elindex=1}^\nel\frac{\rd}{\rd t} ( \|\wf_\elindex-\wf_\elindex' \|^2_{L^2} )\\
&=
2\mathrm{Im} \langle \pot_{\text{Hx}} [\dens ]\wf-\pot_{\text{Hx}} [\dens' ]\wf',\wf-\wf' \rangle_{\ltwoN}\\
&\leq C \|\wf-\wf' \|^2_{\ltwoN},\nonumber
\end{align*} 
where $C=C ( \|\wf \|_{C^0 ([0,\tau];\htwoN)}, \|\wf' \|_{C^0 ([0,\tau];\htwoN)},\tau,q,\lambda,\nel )>0$ stems from the Cauchy--Schwarz inequality and combining \eqref{eq:A} from Lemma \ref{lem:lipschitzconvolution} and \eqref{eq:D} from Lemma \ref{lem:lipschitzlda}. Finally, by Gr\"onwall's lemma we get that $\wf=\wf'$ and this concludes the proof. 
\end{proof}
\begin{lemma}
\label{lem:mappingE}
Let $q\geq7/2$ and $\lambda\in\bb{R}$. Let $\tau>0$ be such that the following holds: for given~$\nuc\in\Bnuc$, $\wf\in\Bel$ is the unique short-time solution to \eqref{eq:KS}. Then, the mapping
\begin{align*}
\mathcal{E}:\nuc\in\mathcal{B}_{\mathrm{nuc}} (\tau )\longmapsto\wf\in\mathcal{B}_{\mathrm{el}} (\tau ),
\end{align*}
is bounded and continuous as a map from $C^0\big([0,\tau];\bb{R}^{3\nnuc}\big)$ to $C^0\big([0,\tau];\ltwoN\big)$.
\end{lemma}
\begin{proof}
Since $\Bel$ is a bounded subset of $C^0 ([0,\tau];\ltwoN )$, the mapping $\mathcal{E}$ is bounded in the $C^0 ([0,\tau];\ltwoN )$ norm. In order to prove continuity of $\mathcal{E}$ as a map from $C^0\big([0,\tau];\bb{R}^{3\nnuc}\big)$ to $C^0\big([0,\tau];\ltwoN\big)$, we consider a sequence $ \{\nuc_n \}_{n\in\bb{N}}\subset\Bnuc$ such that $\nuc_n\xrightarrow{n\longrightarrow\infty}\nuc\in\Bnuc$ in the $C^0 ([0,\tau];\bb{R}^{3\nnuc} )$ norm. Similarly to $\wf=\mathcal{E} [\nuc ]$, we define $\wf_n:=\mathcal{E} [\nuc_n ]$ with $\dens_n:= |\wf_n |^2$. Then,
\begin{align*}
\displaystyle i\frac{\partial}{\partial t} (\wf_n-\wf )=\schop [\nuc,\dens ] (\wf_n-\wf )+{}\zeta{}_n,\qquad
 (\wf_n-\wf )(0)=0,
\end{align*}
with
\begin{align}
{}\zeta{}_n&:={}\zeta{}_n^1+{}\zeta{}_n^2+{}\zeta{}_n^3,\nonumber\\
{}\zeta{}_n^1&:=
\extpot[\nuc_n ]\wf_n-\extpot[\nuc ]\wf-\extpot[\nuc ] (\wf_n-\wf )=
 (\extpot[\nuc_n ]-\extpot[\nuc ] )\wf_n,\nonumber\\
{}\zeta{}_n^2&:=\potH [\dens_n ]\wf_n-\potH [\dens ]\wf-\potH [\dens ] (\wf_n-\wf )\nonumber\\
&=
\sum_{\elindex=1}^\nel
 \{\mathrm{Re} \big[\overline{ (\wf_{n\elindex}-\wf_\elindex )} (\wf_{n\elindex}+\wf_\elindex ) \big]\star  |\cdot |^{-1} \}\wf_n,\label{eq:star18}\\
{}\zeta{}_n^3&:=\potX [\dens_n ]\wf_n-\potX [\dens ]\wf-\potX [\dens ] (\wf_n-\wf )=
\lambda [\dens_n^{q-1}-\dens^{q-1} ]\wf_n,\nonumber
\end{align}
where we used $|a|^2-|b|^2=\mathrm{Re} [\overline{ (a-b )} (a+b ) ]$ in \eqref{eq:star18}. We denote by~$ \{\schop [\nuc(t),\dens ],t\in[0,T] \}$ the family of KS Hamiltonians for the given $\nuc\in\Bnuc$. Note that since $\wf$ and thus $\dens$ are fixed now, these Hamiltonians are acting linearly on~$\wf_n-\wf$, and can thus be written, similarly to \eqref{eq:linHam}, as
\begin{align*}
\schop(t)=-\tfrac{1}{2}\lap + \pot(t)+\pot_{\text{Hx}} [\dens ]
\end{align*}
with $\pot$ from \eqref{eq:calV}. The linear potential $\pot(t)+\pot_{\text{Hx}} [\dens ]$ satisfies Assumption (A.1) of \cite[Theorem 1.1]{Yajima1987ExistenceEquations}; hence, there exists a family of evolution operators $ \{\mathfrak{U}(t,s),(t,s)\in[0,T]^2 \}$, associated with this family of Hamiltonians, satisfying properties (i)---(iv) of our Lemma \ref{lem:evolutionoperators}. By this, from which it follows that for fixed $t\in[0,T]$, $\mathfrak{U}(t,\cdot){}\zeta{}_n\in C^0 ([0,t],\ltwoN )$, and \cite[Cor. 1.2. (4)]{Yajima1987ExistenceEquations}, we can argue like \cite[p. 982]{Cances1999OnDynamics}, and the corresponding integral representation holds for all $t\in[0,T]$:
\begin{align*}
 (\wf_n-\wf )(t)=-i\int_0^t\mathfrak{U}(t,\sigma){}\zeta{}_n(\sigma)\rd \sigma.
\end{align*}
Using Lemma \ref{lem:evolutionoperators} (iii), we bound for all $n\in\bb{N}$ and $t\in[0,\tau]$
\begin{align*}
& \| (\wf_n-\wf )(t) \|_{\ltwoN}\lesssim{}
\sum_{j\in\{1,2,3\}}\int_0^t \|{}\zeta{}_n^j(\sigma) \|_{\ltwoN}\rd\sigma.
\end{align*}
So, now we deduce $\ltwoN$ estimates on ${}\zeta{}_n^j(\sigma)$ for $j\in\{1,2,3\}$ for all $\sigma\in(0,t)$, using that $\wf_n$ and $\wf$ are elements in $\Bel$, which makes them uniformly bounded with respect to $n$ in~$C^0 ([0,\tau]; \htwoN )$. 
For $j=1$, as noted in \cite[p. 982]{Cances1999OnDynamics}, it holds for all $0< \sigma< t\leq \tau \leq T$ that 
\begin{align*}
 \|{}\zeta{}_n^1(\sigma) \|_{\ltwoN}
&\leq  C_{1,n}\xrightarrow{n\longrightarrow\infty}0
\end{align*}
for some $C_{1,n}=C_{1,n} (\alpha,\wf^0 )>0$. For $j=2$, we use the mapping $G$ from \eqref{eq:mappingG}. This gives for all $\sigma\in(0,t)$
\begin{align*}
 \|{}\zeta{}_n^2(\sigma) \|_{\ltwoN}
&\leq 
\sum_{\elindex=1}^\nel
 \|G [\wf_{n\elindex}(\sigma)-\wf_\elindex(\sigma),\wf_{n\elindex}(\sigma)+\wf_\elindex(\sigma) ] \|_{L^\infty} \|\wf_n(\sigma) \|_{\ltwoN}\nonumber\\
&\lesssim{} \sum_{\elindex=1}^\nel
 \|\wf_{n\elindex}(\sigma)-\wf_\elindex(\sigma) \|_{L^2} \|\wf_{n\elindex}(\sigma)+\wf_\elindex(\sigma) \|_{H^2} \|\wf_n \|_{L^\infty ([0,\tau];\ltwoN )}\nonumber\\
&\leq C_2 \|\wf_{n}(\sigma)-\wf(\sigma) \|_{\ltwoN}
\end{align*}
for some $C_2=C_2 (\alpha,\wf^0 )>0$. For $j=3$ and all~$\sigma\in(0,t)$
\begin{align*}
 \|{}\zeta{}_n^3(\sigma) \|_{\ltwoN}
&\stackrel{\eqref{eq:MVE}}{\lesssim}_{q,\lambda} ( \|\dens_n(\sigma) \|_{L^\infty}^{q-3/2}+ \|\dens(\sigma) \|_{L^\infty}^{q-3/2} ) \|\dens_n(\sigma) \|_{L^\infty}^{1/2} \|\wf_n(\sigma)-\wf(\sigma) \|_{\ltwoN}\nonumber\\
&\stackrel{\eqref{eq:rhoinfty}}{\leq} C_3 \|\wf_n(\sigma)-\wf(\sigma) \|_{\ltwoN},
\end{align*}
for some $C_3=C_3 (q,\alpha,\wf^0 )>0$. Combining these three estimates, for all~$t\in [0,\tau]$
\begin{align*}
 \| (\wf_n-\wf )(t) \|_{\ltwoN}\leq C_{1,n}\tau+  (C_2+C_3 )\int_0^t
 \| (\wf_n-\wf)(\sigma) \|_{\ltwoN}\rd\sigma,
\end{align*}
Now, by Gr\"onwall's lemma, we conclude that for all $t\in[0,\tau]$
\begin{align*}
 \| (\wf_n-\wf )(t) \|_{\ltwoN}\leq C_{1,n}\tau e^{ (C_2+C_3 )t},
\end{align*}
and this concludes the proof.
\end{proof}
\section{Proof of Theorem \ref{thm:shorttimeexistence}}
\label{sec:cauchy}
In this section, we prove the main result, Theorem \ref{thm:shorttimeexistence}.
\begin{lemma}
\label{lem:shorttimeexistence}
Let $q\geq 7/2$ and $\lambda\in\bb{R}$. Further, let $\wf^0\in H^2 (\bb{R}^3;\bb{C}^\nel )$, $\vel^0\in\bb{R}^{3\nnuc}$ and $\nuc^0\in\bb{R}^{3\nnuc}$ be given, with $\nuc^0_{\nucindex}\neq \nuc^0_{\nucindextwo}$ for $\nucindex\neq \nucindextwo$.\newline
Then, there exists $\tau>0$ such that the initial-value problem associated with the system \KSCL{} with $\wf(0)=\wf^0$, $\nuc(0)=\nuc^0$ and $\dot{\nuc}(0)=\vel^0$ has a solution $(\wf,\nuc)\in \mathcal{X}(\tau)$.
\end{lemma}
\begin{proof}
Let $\tau>0$ be such that the following statements hold. For given~$\wf\in\Bel$, \eqref{eq:N} has a unique solution $\nuc\in\Bnuc\cap C^2 ([0,\tau];B_\delta (\nuc^0 ) )$, and for given~$\nuc\in\Bnuc$, \eqref{eq:KS} has a unique solution $\wf\in\Bel$. Existence of such $\tau$ has been proven in Lemmas \ref{lem:KStoCL} and \ref{lem:CLtoKS}. We define the inclusion
\begin{align*}
\mathcal{I}:\Bnuc \cap C^2 ([0,\tau];\bb{R}^{3\nnuc} )&\varlonghookrightarrow\Bnuc,
\end{align*}
which is a continuous and compact mapping. Also, we define the mapping
\begin{align*}
\mathcal{K}:\Bnuc\longrightarrow \Bnuc,\qquad \mathcal{K}:=\mathcal{I}\circ \mathcal{N}\circ\mathcal{E}.
\end{align*}
Since by Lemma \ref{lem:KStoCL}, $\mathcal{N}$ is bounded in the $C^1 ([0,\tau];\bb{R}^{3\nnuc} )$ topology, by the Arzel\`a--Ascoli theorem it follows that $\mathcal{K}$ is a compact mapping, where $\Bnuc$ is equipped with the $C^0 ([0,\tau];\bb{R}^{3\nnuc} )$ topology.

By the classical Schauder's fixed point theorem, $\mathcal{K}$ has a fixed point $\nuc$ in $\Bnuc$. Setting $\wf:=\mathcal{E} [\nuc ]$, the corresponding pair $(\wf,\nuc)$ is the desired solution, and this concludes the proof.
\end{proof}
\begin{lemma}
\label{lem:coupledsolution}
Let $q\geq 7/2$ and $\lambda\in\bb{R}$. Let~$(\nuc,\wf), (\nuc',\wf' )\in \mathcal{X}(\tau)$ be two solutions to \KSCL{} for some $\tau>0$. Then, for all $t\in[0,\tau]$
\begin{align}
 | (\ddot{\nuc}-\ddot{\nuc}' )(t) |
&\leq
C [ | (\nuc-\nuc' )(t) |+ \| (\wf-\wf' )(t) \|_{ L^{3,\infty}(\bb{R}^3;\bb{C}^\nel)} ],\label{eq:a}\\
 \| (\wf-\wf' )(t) \|_{ L^{3,\infty}(\bb{R}^3;\bb{C}^\nel)}
&\leq
C\int_0^t\frac{1}{\sqrt{t-\sigma}} [ | (\nuc-\nuc' )(\sigma) |+ \| (\wf-\wf' )(\sigma) \|_{ L^{3,\infty}(\bb{R}^3;\bb{C}^\nel)} ]\rd\sigma,\label{eq:b}
\end{align}
where $C=C ( \|\wf \|_{C^0 ([0,\tau]; \htwoN )}, \|\wf' \|_{C^0 ([0,\tau]; \htwoN )} )$. \end{lemma}
\begin{proof} We focus on justifying each estimate separately, as follows. \\ \noindent {\it Proof of \eqref{eq:a}}.
In this proof, we use for short-hand notation the function $\Xi:\bb{R}^3\longrightarrow \bb{R}^3$ (a.e.),~$\el\longmapsto\el |\el |^{-3}$ again.\\\\For all $t\in[0,\tau]$ and $\nucindex=1,\ldots,\nnuc$
\begin{align*}
 | (\ddot{\nuc}_{\nucindex}-\ddot{\nuc}_{\nucindex}' )(t) |&\leq
 |\acc_{\nucindex}^1 [\dens(t) ] (\nuc(t) )-\acc_{\nucindex}^1 [\dens'(t) ] (\nuc'(t) ) |+
 |\acc_{\nucindex}^2 (\nuc(t) )-\acc_{\nucindex}^2 (\nuc'(t) ) |\\
&\leq 
\mrm{(I)}+\mrm{(II)}+\mrm{(III)},\nonumber\\
\mrm{(I)}&:= |\acc_{\nucindex}^1 [\dens(t) ] (\nuc(t) )
-\acc_{\nucindex}^1 [\dens(t) ] (\nuc'(t) ) |,\nonumber\\
\mrm{(II)}&:= |\acc_{\nucindex}^1 [\dens(t) ] (\nuc'(t) )
-\acc_{\nucindex}^1 [\dens'(t) ] (\nuc'(t) ) |,\nonumber\\
\mrm{(III)}&:= |\acc_{\nucindex}^2 (\nuc(t) )
-\acc_{\nucindex}^2 (\nuc'(t) ) |.\nonumber
\end{align*}
By Lemma \ref{lem:forces} on the force functions, $\acc_{\nucindex}^1 [\dens ]$ are uniformly Lipschitz continuous in the nuclear variable for all $t\in[0,\tau]$ and $\nucindex$, by which
\begin{align*}
\text{(I)}&\lesssim{}\sum_{\elindex=1}^\nel | \langle\wf_\elindex(t),\Xi(\cdot-\nuc_{\nucindex}(t) )\wf_\elindex(t) \rangle_{L^2}
-
 \langle\wf_\elindex(t),\Xi(\cdot-\nuc_{\nucindex}'(t) )\wf_\elindex(t) \rangle_{L^2}
 |\\
&\stackrel{\eqref{eq:forcebounds2}}{\leq} C_{\mrm{I}} | (\nuc_{\nucindex}-\nuc_{\nucindex}' )(t) |\leq C_{\mrm{I}} | (\nuc-\nuc' )(t) |
\end{align*}
for some $C_{\mrm{I}}=C_{\mrm{I}} ( \|\wf \|_{C^0 ([0,\tau]; \htwoN )} )>0$. Also, 
\begin{align*}
\text{(II)}&\leq C_{\mrm{II}}
 \| (\wf-\wf' )(t) \|_{ L^{3,\infty}(\bb{R}^3;\bb{C}^\nel)}
\end{align*}
for some $C_{\mrm{II}}=C_{\mrm{II}} ( \|\wf \|_{C^0 ([0,\tau]; \htwoN )}, \|\wf' \|_{C^0 ([0,\tau]; \htwoN )} )>0$ by H\"older's inequality in Lorentz spaces \eqref{eq:holderlorentz} and the fact that $ \|\cdot \|^{-2}_{\bb{R}^3}\in L^{3/2,\infty}$. Since $\nuc,\nuc'\in\Bnuc$, we bound (III) similarly to Part 1 of the proof of Lemma \ref{lem:KStoCL}:
\begin{align*}
\text{(III)}&\lesssim{}\sum_{\subalign{&\nucindextwo=1,\\&\nucindextwo\neq \nucindex}}^{\nnuc} |\Xi( (\nuc_{\nucindex}-\nuc_{\nucindextwo} )(t) )-\Xi( (\nuc_{\nucindex}'-\nuc_{\nucindextwo}' )(t) ) |\lesssim_{\delta,\nuc^0} | (\nuc-\nuc' )(t) |.
\end{align*}
Since these results hold for all $\nucindex$, \eqref{eq:a} follows.\\\\
\textit{Proof of \eqref{eq:b}.}\\Similarly to the proof of Lemma \ref{lem:CLtoKS},
\begin{align*}
&
\begin{aligned}
\displaystyle i\frac{\partial}{\partial t} (\wf-\wf' )&=-\displaystyle\tfrac{1}{2}\lap(\wf-\wf' )+
\extpot[\nuc ] (\wf-\wf' )
+\pot_{\text{Hx}} [\dens ] (\wf-\wf' )
+\wt{{}\zeta{}},\\
 (\wf-\wf' )(0)&=0,
\end{aligned}
\end{align*}
where $\wt{{}\zeta{}}:=\wt{{}\zeta{}}^1+\wt{{}\zeta{}}^2+\wt{{}\zeta{}}^3$, with for $j\in\{1,2,3\}$, $\wt{{}\zeta{}}^j$ being ${}\zeta{}_n^j$ with $ (\nuc_n,\wf_n )\longmapsto  (\nuc',\wf' )$. As the operator $-\lap/2$ generates the free propagator $U_0$, we write the equivalent integral equation for all $t\in[0,\tau]$
\begin{align*}
 (\wf-\wf' )(t)=-i\int_0^tU_0(t-\sigma) \{\extpot[\nuc(\sigma) ] (\wf-\wf' )(\sigma)
+\pot_{\text{Hx}} [\dens ] (\wf-\wf' )(\sigma)
+\wt{{}\zeta{}}(\sigma) \}\rd\sigma.
\end{align*}
We recall that by \cite[Lemma 6]{Cances1999OnDynamics},  for all $\sigma\in(0,\tau]$ and $f\in L^{3/2,\infty}$, it holds that 
\begin{align*}
 \|U_0(\sigma)f \|_{L^{3,\infty}}\lesssim{} \frac{1}{\sqrt{\sigma}} \|f \|_{L^{3/2,\infty}}.
\end{align*}
Using this estimate, we obtain for all $t\in[0,\tau]$ and $\elindex=1,\ldots,\nel$
\begin{align*}
 \| (\wf_\elindex-\wf_\elindex' )(t) \|_{L^{3,\infty}}&\lesssim{}\int_0^t\frac{1}{\sqrt{t-\sigma}}\Big[ \| (\extpot[\nuc(\sigma) ] (\wf-\wf' ) )_\elindex(\sigma) \|_{L^{3/2,\infty}}\\
&\quad+
 \| (\potH [\dens ] (\wf-\wf' ) )_\elindex(\sigma) \|_{L^{3/2,\infty}}+ \| (\potX [\dens ] (\wf-\wf' ) )_\elindex(\sigma) \|_{L^{3/2,\infty}}
\nonumber\\
&\quad
+\sum_{j\in\{1,2,3\}} \| (\wt{{}\zeta{}}^j(\sigma) )_\elindex \|_{L^{3/2,\infty}}\Big]\rd\sigma.\nonumber
\end{align*}
Since $ \|\cdot \|^{-1}_{\bb{R}^3}\in L^{3,\infty}$, by H\"older's inequality on $L^{3/2,\infty}$, we obtain for all~$\sigma\in(0,t)$ and $\elindex=1,\ldots,\nel$ that
\begin{align*}
 \| (\extpot[\nuc(\sigma) ] (\wf-\wf' ) )_\elindex(\sigma) \|_{L^{3/2,\infty}}&\lesssim{}\sum_{\nucindex=1}^{\nnuc} \| |\cdot-\nuc_{\nucindex}(\sigma) |^{-1} \|_{L^{3,\infty}} \| (\wf_\elindex-\wf_\elindex' )(\sigma) \|_{L^{3,\infty}}\\
&\lesssim{} \| |\cdot |^{-1} \|_{L^{3,\infty}} \| (\wf_\elindex-\wf_\elindex' )(\sigma) \|_{L^{3,\infty}}\nonumber\\
&\lesssim{} \| (\wf_\elindex-\wf_\elindex' )(\sigma) \|_{L^{3,\infty}}.\nonumber
\end{align*}
Note also that
\begin{align}
 \| (\potH [\dens ] (\wf-\wf' ) )_\elindex(\sigma) \|_{L^{3/2,\infty}}&
\lesssim{} \|\dens(\sigma)\star |\cdot |^{-1} \|_{L^{3,\infty}}
 \| (\wf_\elindex-\wf_\elindex' )(\sigma) \|_{L^{3,\infty}}\label{eq:star111}\\
&\lesssim{} \|\dens(\sigma) \|_{L^{1}} \| |\cdot |^{-1} \|_{L^{3,\infty}}
 \| (\wf_\elindex-\wf_\elindex' )(\sigma) \|_{L^{3,\infty}}\label{eq:star211}\\
&\leq C_\text{H}
 \| (\wf_\elindex-\wf_\elindex' )(\sigma) \|_{L^{3,\infty}}\nonumber
\end{align}
for some $C_\text{H}=C_\text{H} ( \|\wf \|_{C^0 ([0,\tau];\htwoN)} )>0$. Here, we have used H\"older's inequality on $L^{3/2,\infty}$ in \eqref{eq:star111} and Young's convolution inequality on $L^{3,\infty}$ in \eqref{eq:star211}. Moreover,
\begin{align}
 \| (\potX [\dens ] (\wf-\wf' ) )_\elindex(\sigma) \|_{L^{3/2,\infty}}
&\lesssim_{\lambda} \| [\dens(\sigma) ]^{q-1} \|_{L^{3,\infty}} \| (\wf_\elindex-\wf_\elindex' )(\sigma) \|_{L^{3,\infty}}\label{eq:star112}\\
&\lesssim{}
 \| [\dens(\sigma) ]^{q-1} \|_{L^{3}} \| (\wf_\elindex-\wf_\elindex' )(\sigma) \|_{L^{3,\infty}}\label{eq:star212}\\
&
\leq C_\text{x}
 \| (\wf_\elindex-\wf_\elindex' )(\sigma) \|_{L^{3,\infty}}\label{eq:star312}
\end{align}
for some $C_\text{x}=C_\text{x} ( \|\wf \|_{C^0 ([0,\tau]; \htwoN )},q )>0$. Here, we used  H\"older's inequality on~$L^{3/2,\infty}$ in \eqref{eq:star112} and \cite[Chapter 4, Prop. 4.2.]{Bennett1988TheTheorems} in \eqref{eq:star212}. In \eqref{eq:star312}, we used Sobolev's inequality with interpolation, and the embedding of $H^2$ into $L^\infty$, by which, with $\theta:=6(q-1)>6$,
\begin{align*}
 \| [\dens(\sigma) ]^{q-1} \|_{L^{3}}^3&\lesssim_{q}
\sum_{\elindex=1}^\nel \| [\wf_\elindex(\sigma) ]^{(\theta-6)+6} \|_{L^1}\\
&\lesssim{}
\sum_{\elindex=1}^\nel \|\wf_\elindex(\sigma) \|_{L^\infty}^{\theta-6} \|\wf_\elindex(\sigma) \|_{L^6}^{6}\lesssim{}
\sum_{\elindex=1}^\nel \|\wf_\elindex \|_{C^0 ([0,\tau];H^2 )}^{\theta}.\nonumber
\end{align*}
On the other hand
\begin{align}
 \| (\wt{{}\zeta{}}^1(\sigma) )_\elindex \|_{L^{3/2,\infty}}
&
\lesssim{}\sum_{\nucindex=1}^{\nnuc}
 \|
 (
 |
\cdot-\nuc_{\nucindex}(\sigma)
 |^{-1}-
 |
\cdot-\nuc_{\nucindex}'(\sigma)
 |^{-1}
 )\wf_\elindex' (\sigma,\cdot ) \|_{L^{3/2,\infty}}\nonumber\\
&
=\sum_{\nucindex=1}^{\nnuc}
 \|
 (
 |
\cdot- (\nuc_{\nucindex}
-\nuc_{\nucindex}' )(\sigma)
 |^{-1}-
 |
\cdot
 |^{-1}
 )\wf_\elindex' (\sigma,\cdot+\nuc_{\nucindex}' ) \|_{L^{3/2,\infty}}\nonumber\\
&\lesssim{}
 \|\wf_\elindex' (\sigma )
 \|_{L^{\infty}}
\sum_{\nucindex=1}^{\nnuc}
 \|
 |
\cdot
 |^{-1}
 |
\cdot- (\nuc_{\nucindex}
-\nuc_{\nucindex}' )(\sigma)
 |^{-1} \|_{L^{3/2,\infty}}\times \nonumber\\
&\qquad\times 
 | (\nuc_{\nucindex}
-\nuc_{\nucindex}' )(\sigma) |\label{eq:star313}\\
&\lesssim{}
 \|\wf_\elindex' \|_{C^0 ([0,\tau];H^2 )}
 \|
 |
\cdot
 |^{-1}
 \|_{L^{3,\infty}}
\sum_{\nucindex=1}^{\nnuc}
 \|
 |
\cdot- (\nuc_{\nucindex}
-\nuc_{\nucindex}' )(\sigma)
 |^{-1}
 \|_{L^{3,\infty}}\times\nonumber
\\
&\qquad\times 
 | (\nuc_{\nucindex}
-\nuc_{\nucindex}' )(\sigma) |\label{eq:star413}\\
&\lesssim_{\nnuc}
 \|\wf_\elindex' \|_{C^0 ([0,\tau];H^2 )}
 \|
 |
\cdot
 |^{-1}
 \|^2_{L^{3,\infty}}
 | (\nuc
-\nuc' )(\sigma) |\nonumber\\
&\lesssim{}
 \|\wf_\elindex' \|_{C^0 ([0,\tau];H^2 )}
 | (\nuc
-\nuc' )(\sigma) |.\nonumber
\end{align}
where we used the triangle inequality written as $ | |\cdot |- |\cdot- (\nuc_{\nucindex}
-\nuc_{\nucindex}' )(\sigma) | |\leq  | (\nuc_{\nucindex}
-\nuc_{\nucindex}' )(\sigma) |$ in \eqref{eq:star313}, H\"older's inequality in $L^{3/2,\infty}$ and the embedding of $H^2$ into~$L^\infty$ in \eqref{eq:star413}. Moreover,
\begin{align}
 \| (\wt{{}\zeta{}}^2(\sigma) )_\elindex \|_{L^{3/2,\infty}}
&\lesssim{}\sum_{\elindextwo=1}^\nel
 \| \{ [\overline{ (\wf_\elindextwo-\wf_\elindextwo' )(\sigma)} (\wf_\elindextwo+\wf_\elindextwo' )(\sigma) ]\star  |\cdot |^{-1} \}\wf_\elindex'(\sigma)
 \|_{L^{3/2,\infty}}\nonumber\\
&\lesssim{}\sum_{\elindextwo=1}^\nel
 \| [\overline{ (\wf_\elindextwo-\wf_\elindextwo' )(\sigma)} (\wf_\elindextwo+\wf_\elindextwo' )(\sigma) ]\star  |\cdot |^{-1}
 \|_{L^{6,\infty}}
 \|\wf_\elindex'(\sigma)
 \|_{L^{2,\infty}}\label{eq:star214}\\
&\lesssim{}\sum_{\elindextwo=1}^\nel
 \| [\overline{ (\wf_\elindextwo-\wf_\elindextwo' )(\sigma)} (\wf_\elindextwo+\wf_\elindextwo' )(\sigma) ]\star  |\cdot |^{-1}
 \|_{L^{6,2}}
 \|\wf_\elindex'(\sigma)
 \|_{L^{2}}\label{eq:star314}\\
&\lesssim{}
 \|\wf_\elindex' \|_{C^0 ([0,\tau];H^{2} )}\sum_{\elindextwo=1}^\nel
 \|
\overline{ (\wf_\elindextwo-\wf_\elindextwo' )(\sigma)} (\wf_\elindextwo+\wf_\elindextwo' )(\sigma)
 \|_{L^{6/5,2}}\times\nonumber\\
&\qquad\times \| |\cdot |^{-1}
 \|_{L^{3,\infty}}\label{eq:star414}\\
&\lesssim{}
 \|\wf_\elindex' \|_{C^0 ([0,\tau];H^{2} )}\times
\nonumber\\
&\qquad \times\sum_{\elindextwo=1}^\nel
 [ \|\wf_\elindextwo(\sigma) \|_{L^{2}}
+ \|\wf_\elindextwo'(\sigma)
 \|_{L^{2}}
 ] \|
 (\wf_\elindextwo-\wf_\elindextwo' )(\sigma)
 \|_{L^{3,\infty}}
\label{eq:star514}\\
&
\leq C_2
 \|
 (\wf-\wf' )(\sigma)
 \|_{ L^{3,\infty}(\bb{R}^3;\bb{C}^\nel)}\nonumber
\end{align}
for some $C_2=C_2 ( \|\wf \|_{C^0 ([0,\tau];\htwoN)},
 \|\wf' \|_{C^0 ([0,\tau];\htwoN)} )>0$.
Here, we used H\"older's inequality on $L^{3/2,\infty}$ in \eqref{eq:star214}, \cite[Chapter 4, Prop. 4.2.]{Bennett1988TheTheorems} in \eqref{eq:star314}, Young's convolution inequality on $L^{6,2}$ in \eqref{eq:star414}, and H\"older's inequality on $L^{6/5,2}$ in \eqref{eq:star514}. Finally,
\begin{align}
 \| (\wt{{}\zeta{}}^3(\sigma) )_\elindex \|_{L^{3/2,\infty}}&\stackrel{\eqref{eq:MVE}}{\lesssim}_{\lambda,q} \|
 [\dens(\sigma) ]^{q-3/2}+ [\dens'(\sigma) ]^{q-3/2}
 \|_{L^{\infty}} \|\wf_\elindex'(\sigma) |\wf(\sigma)-\wf'(\sigma) |
 \|_{L^{3/2,\infty}}\nonumber\\
&\lesssim{}
 \big[ \|\dens(\sigma) \|^{q-3/2}_{L^\infty}+ \|\dens'(\sigma) \|_{L^\infty}^{q-3/2} \big]
 \|\wf_\elindex'(\sigma) \|_{L^{3,\infty}}\times\nonumber \\
&\qquad\times  \| | (\wf-\wf' )(\sigma) |
 \|_{L^{3,\infty}(\bb{R}^3;\bb{C}^\nel)}\label{eq:star215}\\
&\lesssim_{q}
 \Big[
\sum_{\elindexthree=1}^\nel \|\wf_{\elindexthree}(\sigma) \|_{H^2}^{2q-3}
+
\sum_{\elindexfour=1}^\nel \|\wf_{\elindexfour}'(\sigma) \|_{H^2}^{2q-3}
 \Big]\|\wf_\elindex'(\sigma) \|_{L^{3}}
\times\nonumber\\&\qquad\times \sum_{\elindextwo=1}^\nel \| (\wf_\elindextwo-\wf_\elindextwo' )(\sigma)
 \|_{L^{3,\infty}}\label{eq:star315}\\
&\leq C_3 \| (\wf-\wf' )(\sigma)
 \|_{ L^{3,\infty}(\bb{R}^3;\bb{C}^\nel)}.\label{eq:star415}
\end{align}
for some $C_3=C_3 ( \|\wf \|_{C^0 ([0,\tau]; \htwoN )}, \|\wf' \|_{C^0 ([0,\tau]; \htwoN )},q )>0$. Here, we used H\"older's inequality on $L^{3/2,\infty}$ in \eqref{eq:star215}; \eqref{eq:rhoinfty}, \cite[Chapter 4, Prop. 4.2.]{Bennett1988TheTheorems} in \eqref{eq:star315}; and Sobolev's embedding theorem with interpolation in \eqref{eq:star415}. Since all of these estimates hold for all $\sigma\in(0,t)$, and $\elindex=1,\ldots,\nel$, \eqref{eq:b} follows.
\end{proof}
\label{sec:locun}
\begin{proof}[Proof of Theorem \ref{thm:shorttimeexistence}]
Let $\tau>0$ be such that the following statements hold. For given~$\wf\in\Bel$, \eqref{eq:N} has a unique solution $\nuc\in\Bnuc\cap C^2 ([0,\tau];B_\delta (\nuc^0 ) )$, and for given~$\nuc\in\Bnuc$, \eqref{eq:KS} has a unique solution $\wf\in\Bel$. Existence of such $\tau$ has been proven in Lemmas \ref{lem:KStoCL} and \ref{lem:CLtoKS}. Existence of the solution $(\nuc,
\wf)\in \mathcal{X}(\tau)$ to \KSCL{} has been proven in Lemma \ref{lem:shorttimeexistence}. Uniqueness of this solution follows from Lemma \ref{lem:coupledsolution}. For two solutions~$ (\nuc,\wf ), (\nuc',\wf' )\in \mathcal{X}(\tau)$ and $p>2$, let us define the function $h\in C^0 ([0,\tau];\bb{R}_0^+ )$ by
\begin{align*}
h(t):= [ | (\nuc-\nuc' )(t) |+ \| (\wf-\wf' )(t) \|_{ L^{3,\infty}(\bb{R}^3;\bb{C}^\nel)} ]^p.
\end{align*}
Since $\nuc$ and $\nuc'$ both solve \eqref{eq:N} on $[0,\tau]$ and thus are fixed points of the mapping $\mathcal{T}$ in \eqref{eq:mappingT}, for all $t\in[0,\tau]$
\begin{align*}
 | (\nuc-\nuc' )(t) |&\leq
\int_0^t(t-\sigma) | (\ddot{\nuc}-\ddot{\nuc}' )(\sigma) |\rd\sigma.
\end{align*}
Now, using this in combination with Lemma \ref{lem:coupledsolution} in \eqref{eq:star116} and H\"older's inequality, for all $t\in[0,\tau]$
\begin{align}
h(t)&\lesssim_{p}C\biggr\{\int_0^t \Big(t-\sigma+\frac{1}{\sqrt{t-\sigma}} \Big)\times \nonumber\\
&\qquad\times [ | (\nuc-\nuc' )(\sigma) |+ \| (\wf-\wf' )(\sigma) \|_{ L^{3,\infty}(\bb{R}^3;\bb{C}^\nel)} ]\rd\sigma\biggr\}^p\label{eq:star116}\\
&\lesssim{}
C \Big\| \Big(t-\cdot + \frac{1}{\sqrt{t-\cdot}} \Big)h^{1/p} \Big\|^p_{L^1 ([0,t];\bb{R} )}
\nonumber\\
&\lesssim{}
C \Big\|t-\cdot + \frac{1}{\sqrt{t-\cdot}} \Big\|^p_{L^{p'} ([0,t];\bb{R} )}
 \|h^{1/p} \|^p_{L^p ([0,t];\bb{R} )}\lesssim_{\tau}
C\int_0^t h(\sigma)\rd\sigma,\nonumber
\end{align}
where $C=C ( \|\wf \|_{C^0 ([0,\tau]; \htwoN )}, \|\wf' \|_{C^0 ([0,\tau]; \htwoN )}
 )$ is from Lemma \ref{lem:coupledsolution}. Now, using Gr\"onwall's inequality, we obtain $h\leq 0$ on $[0,\tau]$. Since $h\geq 0$ too by definition, and~$h(0)=0$ since $\nuc(0)=\nuc'(0)=\nuc^0$ and $\wf(0)=\wf'(0)=\wf^0$, we get $h\equiv 0$, by which~$ (\nuc,\wf )= (\nuc',\wf' )$. This completes the proof. 
\end{proof}
\section{Acknowledgements}
B.B., O.\c{C}. and W.S. acknowledge funding by the Innovational Research Incentives Scheme Vidi of the Netherlands Organisation for Scientific Research (NWO) with project number 723.016.002.
C.M. was partially funded by The CHERISH Digital Economy Centre
- Collaboration and Knowledge Exchange Support ref 69M. C.M. would like to thank the members of the Department of Mathematics and Computer Science and of the Institute for Complex Molecular Systems - Eindhoven University of Technology, for the warm hospitality.
\section{CRediT author statement}
Wouter Scharpach: \textit{Conceptualisation, formal analysis, writing original draft, review \&{} editing, project administration}; Carlo Mercuri: \textit{Conceptualisation, methodology, formal analysis, writing original draft, review, supervision, funding acquisition}; Bj\"orn Baumeier: \textit{Conceptualisation, review, supervision, funding acquisition}; Mark Peletier: \textit{Conceptualisation, methodology, review}; Georg Prokert: \textit{Conceptualisation, review (in part)}; Onur \c{C}aylak: \textit{Conceptualisation}.
\appendix
\section{Notation}
\label{sec:notation}
Throughout the paper, we make use of the following notation:
\begin{itemize}
    \item We use $A\lesssim{} B$ to denote $ |A |\leq CB$, where $0<C<\infty$.
\item We use $A\lesssim_{\alpha,\beta} B$ to denote dependence on parameters $\alpha,\beta$: so, $ |A |\leq C_{\alpha,\beta}B$ with $0<C_{\alpha,\beta}<\infty$. In our notation, $A\lesssim_{\alpha} B\lesssim_{\beta} \Gamma$ means that $A\lesssim_{\alpha,\beta} \Gamma$.
\end{itemize}
Further, we make use of the following normed spaces:
\begin{itemize}
\item We use the notation $L^{p,r}(\bb{R}^3)=L^{p,r}$, $p\in[1,\infty)$, $r\in[1,\infty]$, for Lorentz spaces, with
\begin{itemize}
\item the radial non-increasing rearrangement of measurable functions $f$ on $\bb{R}^3$
\begin{align*}
f^*(t)=\inf \{s>0 \rvert  | \{\el\in\bb{R}^3 | |f(\el) |>s \} |\leq t \},
\end{align*}
\item the average of $f^*$
\begin{align*}
f^{**}(t)=\frac{1}{t}\int_0^t f^*(s)\rd s,
\end{align*}
and
\item the norms
\begin{align*}
 \|f \|_{L^{p,r}}^r=\int_0^\infty |t^{1/p}f^{**}(t) |^r \frac{\mrm{d}t}{t}
\end{align*}
on $L^{p,r}$, $r<\infty$, and
\begin{align*}
 \|f \|_{L^{p,\infty}}=\sup_{t\in\bb{R}}  |t^{1/p}f^{**}(t) | 
\end{align*}
on the weak Lebesgue spaces $L^{p,\infty}$.
\end{itemize}
\item We use the notation $W^{k,p}(\bb{R}^3)=W^{k,p}$ and $W^{k,2}=H^k$, $k\in\bb{N}$, $p\in[1,\infty]$, for classical Sobolev spaces, setting in particular
\begin{align*}
 \|f \|_{H^2}^2= \|f \|^2_{L^2}+ \|\Delta f \|^2_{L^2}.
\end{align*}
\end{itemize}
Further, we make use of the following inequalities:
\begin{itemize}
\item H\"older's inequality on Lorentz spaces. \cite{ONeil1963ConvolutionSpaces} Let $f\in L^{p_1,q_1}$, $g\in L^{p_2,q_2}$, with $p_1,p_2\in(0,\infty)$, $q_1,q_2\in(0,\infty]$. Then 
\begin{align}
\|f\cdot g\|_{L^{r,s}}\lesssim_{p_1,p_2,q_1,q_2}\|f\|_{L^{p_1,q_1}}\|g\|_{L^{p_2,q_2}}
\label{eq:holderlorentz}
\end{align}
with $1/r=1/p_1+1/p_2$, $1/s=1/q_1+1/q_2$.
\item Young's convolution inequality on Lorentz spaces. \cite[Thm. 2.10.1]{Ziemer1989WeaklyVariation} Let $f\in L^{p_1,q_1},g\in L^{p_2,q_2}$, with $1/p_1+1/p_2>1$. Then 
\begin{align*}
 \|f\star g \|_{L^{r,s}}\leq 3r
 \|f \|_{L^{p_1,q_1}} \|g \|_{L^{p_2,q_2}}
\end{align*}
with $1/r=1/p_1+1/p_2-1$ and $s\in[1,\infty]$ such that $1/q_1+1/q_2\geq 1/s$.
\item Hardy's inequality: \cite{Hardy1952Inequalities}
\begin{align}
 \rVert  |\el-\cdot |^{-1}f \rVert_{L^2}\leq 2 \rVert\nabla f \rVert_{L^2 }\label{eq:hardy}
\end{align}
for all $f\in H^1$ and $\el\in\bb{R}^3$.
\end{itemize}
\bibliographystyle{plain}

\end{document}